\documentclass{nMCM2e}

\usepackage{subfigure}
\usepackage{graphicx}
\usepackage{amsmath,amssymb}
\usepackage{booktabs}
\usepackage{color}
\usepackage{txfonts}

\newcommand{\e}{\mathrm{e}}
\newcommand{\f}{\mathrm{f}}

\theoremstyle{plain}

\theoremstyle{definition}

\theoremstyle{remark}

\begin{document}



\title{{\itshape On fundamental unifying concepts for trajectory--based slow invariant
  attracting manifold computation in multiscale models of chemical kinetics}}

\author{D. Lebiedz$^{\rm a}$$^{\ast}$\thanks{$^\ast$Corresponding author. Email: dirk.lebiedz@uni-ulm.de
\vspace{6pt}} and J. Unger$^{\rm a}$\\\vspace{6pt}  $^{a}${\em{Institute for Numerical Mathematics, 
	Ulm University, Germany}}}

\maketitle

\begin{abstract}
Chemical kinetic models in terms of ordinary differential equations
 correspond to finite dimensional dissipative 
 dynamical systems involving a multiple time scale structure. Most dimension
 reduction approaches aimed at a 
 slow mode--description of the full system compute approximations of
 low--dimensional attracting slow invariant manifolds and parameterize
 these manifolds in terms of a subset of chosen chemical species, the reaction progress variables. The invariance
 property suggests a slow invariant manifold to be constructed as
 (a bundle of) solution trajectories of suitable ordinary differential equation initial or boundary value problems.  
 The focus of this work is on a discussion of fundamental and unifying geometric and analytical issues
 of various approaches to trajectory-based numerical approximation techniques
 of slow
 invariant manifolds that are in practical use for model reduction in chemical kinetics.
 Two basic concepts are pointed out reducing various model reduction
 approaches to a common denominator. In particular, we discuss our recent
 trajectory optimization approach in the light of these two concepts. We relate
 both of them in a variational boundary value viewpoint, propose a Hamiltonian
 formulation and conjecture its relation to conservation laws, (partial)
 integrability and symmetry issues as underlying fundamental principles and
 potentially unifying elements of diverse dimension reduction approaches. 


\begin{classcode}34B15, 34D15, 34C45, 37C45, 37J15, 49S05\end{classcode}
\end{abstract}

\section{Introduction}

A large number of technical applications include chemically reacting
flows comprising an interplay between convective and diffusive species
transport and chemical reaction processes. 
Due to the large number of chemical species involved and the stiffness of the
kinetic ordinary differential equations (ODE) with time scales ranging from
nanoseconds to seconds, simulation of 
chemically reacting flows (for instance in combustion processes) is often nearly
impossible in reasonable computing time, whenever using detailed reaction
mechanism models. 
This calls for multiscale approaches and appropriate model reduction
techniques. 

In chemical reaction kinetics modeled by dissipative ODE systems with spectral
gaps it is observed that solution trajectories
bundle near invariant manifolds of successively lower
dimension during time evolution caused by the multiple time scales
generating spectral gaps. This time scale separation of the model solution
into fast and slow modes 
is the basis for most model and complexity reduction techniques, where the
long time scale system dynamics is approximated via elimination of the fast
relaxing modes by enslaving them to the slow ones. 
The outcome of this is in the ideal case an invariant manifold of slow motion
(denoted as slow invariant manifold (SIM)) possessing the property of
attracting system trajectories from arbitrary initial values. 
Many model reduction methods make use of a species reconstruction
technique for SIM approximation which is provided by an implicitly defined
function mapping a subset of the chemical species of the full model---called reaction progress variables (RPVs)---onto the full species composition by
determining a point on a SIM. 

Among the first model reduction methods have been the Quasi--Steady--State
Assumption (QSSA) \cite{Bodenstein1913,Chapman1913} and the Partial
Equilibrium Approximation (PEA) \cite{Michaelis1913}. In the QSSA approach
certain species are supposed to be in steady state whereas in the PEA approach
several (fast) reactions are assumed to be in equilibrium. 
Due to their conceptual simplicity both the QSSA and the PEA
are still used nowadays although more sophisticated model reduction methods have been developed.
In 1992 Maas and Pope introduced the Intrinsic Low Dimensional Manifold (ILDM) method \cite{Maas1992} which has become
very popular and widely used in the reactive flow community, in particular in
combustion applications. However, both QSSA/PEA manifolds and ILDMs are not
invariant. Other popular techniques are Computational Singular Perturbation
(CSP) \cite{Lam1985,Lam1994} proposed by Lam in 1985, the method described in \cite{Kevrekidis2003,Theodoropoulos2000}, 
which is based on equation-free approaches, and in \cite{Chiavazzo2011}, the relaxation redistribution method has recently been published.
Mease et al. describe a method, where finite-time Lyapunov exponents and
vectors are used to determine slow manifolds \cite{Mease2012}.
The use of those manifolds for the numerical solution of differential equations is obvious in the G-scheme \cite{Valorani2009}.
Furthermore, there are other approaches, whereof a few are presented in Chapter \ref{chapterthree} of this work, including the 
Invariant Constrained Equilibrium Edge PreImage Curve (ICE-PIC) introduced by
Ren et al.\ in 2006 \cite{Ren2005,Ren2006a}, an iterative model reduction method
called Zero--Derivative Principle (ZDP) presented in \cite{Gear2005,Zagaris2009}, the Flow Curvature Method proposed by Ginoux in \cite{Ginoux2008}, 
the Functional Equation Truncation (FET) approach by Roussel \cite{Roussel2006,Roussel2012}, and
methods by Adrover et al.\ \cite{Adrover2007a,Adrover2007} and Al-Khateeb et al.\ \cite{Al-Khateeb2009}.
Additionally, there is a trajectory--based optimization approach based on a
variational principle proposed by Lebiedz et al., described in
\cite{Lebiedz2004c,Lebiedz2006b,Lebiedz2009,Lebiedz2010a,Lebiedz2011,Lebiedz2011a,Reinhardt2008}.

The aim of the present work is the discussion of common ideas and concepts underlying diverse
kinetic model reduction approaches. In particular, various points of view on the
latter variational principle are presented
yielding a deeper insight into its fundamental ideas and its relation to
other techniques for SIM identification.

\section{Slow invariant manifolds for 2-D test models: analytical treatment} \label{anatrea}

This paper considers two ODE test models for SIM computation, both well
suited for analytical treatment due to the availability of explicit formula
for the SIM to be approximated by various model reduction techniques.
Those systems are two--dimensional allowing a clear visualization of phase
space dynamics. 
The first system is the linear model 
\begin{subequations}\label{linsys}
\begin{align}
 \partial_tz_1(t)&=\left(-1-\frac{\gamma}{2}\right)z_1(t)+\frac{\gamma}{2}z_2(t)\\
 \partial_tz_2(t)&=\frac{\gamma}{2}z_1(t)+\left(-1-\frac{\gamma}{2}\right)z_2(t),\quad\gamma>0,
\end{align}
\end{subequations}
with $\gamma\in\mathbb{R}$, $t\in\mathbb{R}$, and $z_1,z_2\in C^{\infty}\left(\mathbb{R},\mathbb{R}\right)$, whereas the second system is a nonlinear one, the Davis--Skodje model
\cite{Davis1999,Singh2002} 
\begin{subequations} \label{DS}
 \begin{align}
  \partial_t z_1(t) &= -z_1(t)\\
  \partial_t z_2(t) &= -\gamma z_2(t)+\frac{(\gamma-1)z_1(t)+\gamma\left(z_1(t)\right)^2}{\left(1+z_1(t)\right)^2},\quad \gamma>1.
 \end{align}
\end{subequations}
In both cases $\gamma$ measures the spectral gap in the system
meaning that varying this parameter modifies the degree of attraction of the
SIM. In the linear system the SIM obviously corresponds to the slow
eigenspace of the system matrix, for the Davis--Skodje model the SIM is the
stable manifold of the fixed point $(z_1,z_2)=(0,0)$ which is tangent to the
slow eigenspace of the system Jacobian evaluated at $(z_1,z_2)=(0,0)$. 
General analytical solutions are given by
\begin{subequations}\label{ansol}
 \begin{align}
  z_1(t) &= c_1\e^{-t}+c_2\e^{(-1-\gamma) t}\\
  z_2(t) &= c_1\e^{-t}-c_2\e^{(-1-\gamma) t},\quad c_1,c_2\in\mathbb{R}\label{ansolnrpv}
 \end{align}
\end{subequations}
for the linear model and by
\begin{subequations}\label{ansolDS}
 \begin{align}
  z_1(t) &= c_1\e^{-t}\label{1}\\
  z_2(t) &= c_2\e^{-\gamma t}+\frac{c_1}{c_1+\e^t},\quad c_1,c_2\in\mathbb{R}\label{4}
 \end{align}
\end{subequations}
for the nonlinear Davis--Skodje model, with integration constants $c_i,\ i=1,2,$ to be
determined by setting initial values. 
For the linear system (\ref{linsys}), trajectory--based model reduction in
terms of eliminating the fast modes implies here setting $c_2$ equal to zero
which leads to    
\begin{subequations}\label{redansol}
 \begin{align}
  z_1(t) &= c_1\e^{-t}\\
  z_2(t) &= c_1\e^{-t}\quad \label{redansolnrpv}
 \end{align}
\end{subequations}
resulting in the SIM $z_1\equiv z_2$.
The same procedure applied to the Davis--Skodje model yields
$z_2(t)=\tfrac{z_1(t)}{z_1(t)+1}$ for the SIM.

In general, an autonomous kinetic ODE model can be formulated as 
\begin{equation}
\partial_tz(t)=S\left(z(t)\right),\quad S\in C^{\infty}\left(\mathbb{R}^n,\mathbb{R}^n\right)
\end{equation}
with $z(t)=\left(z_i(t)\right)_{i=1}^n\in C^{\infty}\left(\mathbb{R},\mathbb{R}^n\right)$ modeling the full
composition state vector of the system. 
Most kinetic model reduction approaches define a subset of the state
variables, the RPVs 
\begin{equation}
 z_j(t),\quad j\in I_{\text{fixed}},
\end{equation}
where $I_{\text{fixed}}\subset\{1,\dots,n\}$ is the index set for the RPVs that parameterize
a SIM. Species reconstruction is the process of fixing the RPVs at a given point in time $t=t_*$ and determining the free variables (non RPVs) 
$z_j(t_*),\ j\notin I_{\text{fixed}}$ from $z_j^{t_*}\coloneqq z_j(t_*),\ j\in
I_{\text{fixed}}$, which implicitly defines a species reconstruction function
$h\in C^{\infty}(\mathbb{R}^{\#I_{\text{fixed}}},\mathbb{R}^{n-\#I_{\text{fixed}}})$ mapping the RPVs to the full species composition thus
determining a point on a SIM. 

Regarding the linear model (\ref{linsys}) and its solution (\ref{ansol}) with $I_{\text{fixed}}=\{2\}$, elimination of the fast mode leads to
$c_1=z_2^{t_*}\e^{t_*}$, that is 
\begin{align}
 z(t_*)=\begin{pmatrix}
          z_1(t_*)=h\left(z_2^{t_*}\right)=z_2^{t_*}\\\\z_2^{t_*}
         \end{pmatrix}
\end{align}
with $h(\cdot)$ being the species reconstruction function. We
call $z(t_*)$ the point of interest (POI) in the following.

Model reduction by species reconstruction for initial value problem (\ref{linsys}) together with $z_1(t_*)=z_1^{t_*},\ z_2(t_*)=z_2^{t_*}$ yields 
\begin{subequations}\label{reducedmodel}
\begin{align}
 \partial_tz_2(t)&=\frac{\gamma}{2}h\left(z_2(t)\right)+\left(-1-\frac{\gamma}{2}\right)z_2(t)=-z_2(t)\label{redlinmod}\\
 z_2(t_*)&=z_2^{t_*}\\
 z_1(t)&=h\left(z_2(t)\right)=z_2(t), \label{additional}
\end{align}
\end{subequations}
with analytical solution of the reduced model equation (\ref{redlinmod}) 
\begin{subequations}
 \begin{align}
  z_2(t) &= c_1\e^{-t}\\
  z_1(t) &= c_1\e^{-t}.
 \end{align}
\end{subequations}
It is obious that it coincides with the general analytical solution of the full
model (\ref{ansol}) where the fast modes $c_2\e^{(-1-\gamma)t}$ have been
eliminated a priori and the initial value has been chosen on the manifold defined by
the function $h(\cdot)$. Thus, the model (\ref{linsys}) together with $z_1(t_*)=z_1^{t_*},\ z_2(t_*)=z_2^{t_*}$ corresponds to the reduced system (\ref{reducedmodel}) for $z_1^{t_*}=h\left(z_2^{t_*}\right)=z_2^{t_*}$.

The same arguments hold for the nonlinear Davis--Skodje test problem (\ref{DS}), (\ref{ansolDS}).

\section{Two basic concepts for slow manifold computation} \label{chapterthree}

Finding a functional with $\Phi(z)=0,\ \Phi\in
C^\infty\left(C^\infty\left(\mathbb{R},\mathbb{R}^n\right),\mathbb{R}^p\right),\ p\leq n$ 
that automatically eliminates the fast modes
(in the case of the two test examples from above a criterion that yields
$c_2=0$) without knowing the analytical solution $z \in
C^\infty\left(\mathbb{R},\mathbb{R}^n\right)$ of the
underlying ODE model equations is the main challenge of trajectory--based model reduction approaches. 
The resulting general species reconstruction problem can be formulated as
\begin{subequations}
 \begin{align}
  \Phi(z)&=0\\
  \partial_tz(t)&=S\left(z(t)\right)\label{zdp_dyn}\\
0 &= g\left(z(t_{*})\right)\label{zdp_ac}\\
  z_j(t_*) &= z_j^{t_*},\quad j\in I_{\text{fixed}}\label{zdp_rpv},
 \end{align}
\end{subequations}
with (\ref{zdp_dyn}) describing the kinetic model equations and (\ref{zdp_rpv}) the fixation of the RPVs at time $t=t_*$.
The function $g\in C^{\infty}\left(\mathbb{R}^n,\mathbb{R}^b\right)$ in (\ref{zdp_ac}) contains possible additional constraints
(for instance chemical element mass conservation relations) and can be omitted for the
two test models above. 
As mentioned in the introduction, there are plenty of different approaches to
find a criterion $\Phi(z)=0$ that (approximately) eliminates the fast modes of the system to obtain a POI that identifies the slow modes and thus the SIM representing the reduced system. 
Since all methods share the same objective, there should be some basic concepts underlying, combining, and collecting different approaches.
The focus of this chapter is the discussion of such basic concepts.

\subsection{Derivative of the state vector}

The first concept various model reduction approaches make use of time
derivatives of the state vector, i.e. $\Phi(z)$ containing terms of type
\begin{align}
 \Phi(z)=\partial^m_tz(t),\quad m\in\mathbb{N},\quad m\geq1.
\end{align}
In the following, some of these methods are briefly reviewed to
  be able to place them in a context.

\subsubsection{Zero--Derivative Principle (ZDP)}\label{ZDP}

A particular species reconstruction method is discussed in 
\cite{Gear2005,Zagaris2009}, which annulates the derivatives of the non RPVs motivating the name
Zero--Derivative Principle (ZDP). The ZDP can in our above framework be fomulated as 
\begin{subequations} \label{zdp}
 \begin{align} \Phi(z)\coloneqq \partial_t^mz_j(t)\Big|_{t=t_*}&=0,\quad j\notin I_{\text{fixed}} \label{zdp_of}\\
    \partial_t z(t) &= S\left(z(t)\right) \\
   0 &= g\left(z(t_{*})\right) \\
   z_j(t_{*}) &= z_j^{t_*},\quad j \in I_{\text{fixed}}.  
 \end{align}
\end{subequations}
The POI $z(t_*)$ is intended to lie in a small neighborhood of a
SIM which is approached with increasing derivative--order $m$ and reached in
the limit $m\rightarrow\infty$ as it is shown in \cite{Gear2005,Zagaris2009}. 

By the help of the linear model (\ref{linsys}), the operation of the ZDP can
be illustrated. Remember, that
elimination of the fast modes characterizes the SIM. Considering the general analytical solution of
(\ref{linsys}) 
\begin{subequations}\label{anasol}
 \begin{align}
  z_1(t) &= c_1\e^{-t}+c_2\e^{(-1-\gamma) t},\quad c_1,c_2\in\mathbb{R}\\
  z_2(t) &= c_1\e^{-t}-c_2\e^{(-1-\gamma) t},\label{anasolnrpv}
 \end{align}
\end{subequations}
the fast modes of the system are represented by the second term of the sum $c_2\e^{(-1-\gamma) t}$ since it includes
the `fast eigenvalue' $-1-\gamma$. Besides the fixation of the RPV, an additional constraint $\Phi\left(z(t)\right)=0$ is needed for obtaining a specified
trajectory---ideally leading to the
elimination of the fast modes ($c_2=0$). 
This is achieved via the zero of the $m^{\text{th}}$ derivative of $z_1$
\begin{align}
 \partial_t^mz_1(t)=(-1)^mc_1\e^{-t}+(-1-\gamma)^mc_2\e^{(-1-\gamma) t}
\end{align}
(in the limit $m\rightarrow\infty$)
where the corresponding eigenvalues of each mode is taken to the power of $m$.
Solving the equation 
\begin{align}
\partial_t^mz_1(t)\Big|_{t=t_*}=0
\end{align}
for $c_2$ ($c_1$ is fixed by the fixation of the RPV) yields
\begin{align}
 c_2=\frac{(-1)^m}{(-1-\gamma)^m}\cdot \text{R}
\end{align}
with R being independent of $m$.
Since $|-1-\gamma|>|-1|$, $c_2\to 0$ for $m\to\infty$, meaning a decreasing contribution of the fast mode for increasing value of $m$.
The same arguments hold for the nonlinear Davis--Skodje test model (\ref{DS})
from \cite{Davis1999,Singh2002}. 
This demonstrates fast modes elimination via the ZDP approach.

\subsubsection{Flow Curvature Method (FCM)}

Another model reduction method based on
derivative--of--the--state--vector--concept is the Flow Curvature Method
proposed by Ginoux \cite{Ginoux2008}, comprising a species reconstruction
technique that computes a special (n-1)-dimensional manifold, the flow curvature manifold, which is defined by the location of the points at which the flow curvature vanishes.
For an $n$-dimensional dynamical system, the zero point of the flow curvature of a trajectory curve is defined as
\begin{align}
 \det{(\partial_tz(t), \partial_t^2z(t),\partial_t^3z(t),\dots,\partial_t^nz(t))}=0 \label{FCM}
\end{align}
with $z(t)\in\mathbb{R}^n$, $n\in\mathbb{N}$. 
Replacing the flow curvature by its successive Lie derivatives (in a
two-dimensional system it is defined by the determinant of the first and third
time derivatives) yields successively higher order approximations of the
SIM (see \cite{Rossetto1986}).
For singularly--perturbed systems, i.e.\ systems comprising a small parameter
$0<\varepsilon<<1,\ \varepsilon\in\mathbb{R}$ controlling the time scale
separation of the system, the analytical equation of the slow manifold
resulting from matched asymptotic expansion in Singular Perturbation Theory \cite{Jones1994, Kaper1999}, which is given by a regular perturbation expansion in $\varepsilon$, equates with the FCM up to a suitable order in $\varepsilon$.
The invariance property of the flow curvature manifold can be shown via the Darboux~Invariance~Theorem~\cite{Darboux1878}.

\subsubsection{Intrinsic Low Dimensional Manifold (ILDM)/Functional Equation
  Truncation (FET)}

As mentioned in the introduction, Maas and Pope introduced the widely used
IDLM method in 1992 \cite{Maas1992} where a local time scale analysis is
performed via matrix decomposition of the Jacobian of the right hand side of the ODE system.
In \cite{Roussel2006,Roussel2012}, Roussel could demonstrate the coincidence
between the ILDM method and his FET approach. The operation concept of FET is shown for a planar system
\begin{subequations}\label{FET}
\begin{align}
 \partial_tz_1(t)&=S_1(z_1(t),z_2(t))\\
 \partial_tz_2(t)&=S_2(z_1(t),z_2(t)).
\end{align}
\end{subequations}
The functional equation
\begin{align}
 S_2(z_1(t),z_2(t))=z_2'(t)S_1(z_1(t),z_2(t)) \label{FETFE}
\end{align}
is achieved by substituting $\partial_tz_2(t)=z_2'(t)S_1(z_1(t),z_2(t))$ with $z_2'(t)=\tfrac{\d z_2(t)}{\d z_1(t)}$ into (\ref{FET}b).
Differentiation of the functional equation (\ref{FETFE}) with respect to $z_1(t)$ yields
\begin{align}\label{FET2}
 S_2'(z_1(t),z_2(t))=z_2''(t)S_1(z_1(t),z_2(t))+z_2'(t)S_1'(z_1(t),z_2(t)).
\end{align}
Motivated by the observation that the error in the ILDM method is
  directly related to the neglect of curvature \cite{Kaper2002}, which is proportional to $z_2''(t)$ here, equation
  (\ref{FET2}) becomes \begin{align}
 z_2'(t)\left(\partial_{z_1}S_1(z_1(t),z_2(t))+z_2'(t)\partial_{z_2}S_1(z_1(t),z_2(t))\right)=\partial_{z_1}S_2(z_1(t),z_2(t))+z_2'(t)\partial_{z_2}S_2(z_1(t),z_2(t)) \label{FETTE}
\end{align}
which is called the truncated equation. 
Thus, we now have two equations ((\ref{FETFE}) and (\ref{FETTE})) in the two unknowns $z_2(t),z_2'(t)$ for every $z_1(t)$
allowing the computation of an approximation to the one-dimensional manifold by using an iterative method to solve (\ref{FETFE}).
The resulting manifold is called Functional Equation Truncation
Approximated (FETA) manifold \cite{Roussel2006,Roussel2012}. Its approximation of the one-dimensional SIM is valid insofar as $z_2''(t)$ is small. 

Already Kaper and Kaper \cite{Kaper2002} pointed out the direct relation
between the ILDM method and the neglect of the curvature, which is used as a
central idea in the above FET approach. The concept of the zero point of a
derivative of the state vector is used here as well.

\subsubsection{Stretching--Based Diagnostics}\label{adrover}

Adrover et al.\ \cite{Adrover2007a,Adrover2007} presented a method for model reduction which
is based on a geometric characterization of local tangent and normal dynamics. 
This description finds its justification in the fact that
the flow along a slow manifold is
slower than the attraction/repulsion to/from it. 
The method uses a ratio $r>1$ of the local stretching (contraction)
rates of vectors orthogonal to the SIM compared to those tangent to the
SIM. Then this ratio is maximized w.r.t.\ $z$. 
Again a two-dimensional dynamical system is considered for demonstration
\begin{equation*}
 \partial_tz(t)=S(z(t))=\begin{pmatrix}S_1(z(t))\\S_2(z(t))\end{pmatrix},\quad z(t)\in\mathbb{R}^2
\end{equation*}
possessing a one-dimensional SIM $\mathcal{W}$.
Then the stretching ratio $r$ is given by
\begin{equation*}
 r(z(t))\coloneqq\frac{\omega_\nu(z(t))}{\omega_\tau(z(t))}\coloneqq\frac{\langle J_S(z(t))\cdot \hat{n}(z(t)),\hat{n}(z(t))\rangle}
{\langle J_S(z(t))\cdot \hat{S}(z(t)),\hat{S}(z(t))\rangle},\quad z(t)\in\mathcal{W}
\end{equation*}
with $\hat{S}(z(t))\coloneqq\tfrac{S(z(t))}{\lVert S(z(t))\rVert}$, $\hat{n}(z(t))\coloneqq\tfrac{n(z(t))}{\lVert n(z(t))\rVert}$, $n(z(t))\coloneqq(S_2(z(t)),-S_1(z(t)))^\top$,
$\langle\cdot,\cdot\rangle$ being the scalar product, $\lVert\cdot\rVert$
indicating the Euclidean norm, and $J_S(z(t))$ being the
Jacobian of the right hand side $S(z(t))$ evaluated at $z(t)$.
Here, $\omega_\tau(z(t))$ denotes the tangential stretching rate and $\omega_\nu(z(t))$ the normal stretching rate.
The reduction method can be viewed as a local embedding technique: Locally
projecting the dynamics onto the slowest directions.
In the $n-$dimensional case ($n>2$) the tangential stretching rate is still given by
\begin{equation*}
 \omega_\tau(z(t))=\langle J_S(z(t))\cdot \hat{S}(z(t)),\hat{S}(z(t))\rangle
\end{equation*}
while the definition of normal stretching rates is
\begin{equation*}
 \omega_\nu(z(t))=\max\limits_{\hat{n}\in N\mathcal{W}_z, \|\hat{n}\|=1} \langle J_S(z(t))\cdot \hat{n}(z(t)),\hat{n}(z(t))\rangle
\end{equation*}
where the maximum is taken over all vectors $\hat{n}(z(t))$ belonging to the normal space $N\mathcal{W}_z$ at $z(t)$.
This value can be computed by the largest eigenvalue of a symmetric
matrix (cf.\ \cite{Adrover2007}).

Since $\partial_t^2z(t)=J_S(z(t))\cdot S(z(t))$ holds, obviously a derivative--of--the--state--vector--concept is used in this method.

\subsubsection{Quasi Steady State Assumption (QSSA)/Partial Equilibrium Approximation (PEA)}

Even in the simplest model reduction approaches QSSA \cite{Bodenstein1913, Chapman1913} and PEA \cite{Michaelis1913} the idea of taking the zero point of a derivative of species can be found.
In the QSSA approach certain species are supposed to be in steady state,
meaning the zero point of the first time derivative of these certain species is regarded.
The correlation between the two model reduction methods is analyzed in
\cite{Goussis2012}, where it is shown that QSSA can be interpreted as a limiting case of PEA.

\subsubsection{Trajectory Based Optimization Approach} \label{tboa}


In \cite{Lebiedz2004c} and follow-up publications of Lebiedz et al.\ \cite{Lebiedz2006b,Lebiedz2009,Lebiedz2010a,Lebiedz2011,Lebiedz2011a,Reinhardt2008} a species
reconstruction method for identifying SIMs is presented based on a variational
principle exploiting trajectory--based optimization
\begin{subequations} \label{op}
 \begin{equation}  \min_{z(t)} \int_{t_{0}}^{t_{\rm{f}}}\!\Phi\left(z(t)\right) \; \textrm{d} t,\quad t_0,t_\f\in\mathbb{R},\quad t_0<t_\f \label{op_of}
 \end{equation}
 \textrm{subject to}
 \begin{align}
   \partial_t z(t) &= S\left(z(t)\right) \label{op_dyn}\\
   0 &= g\left(z(t_{*})\right) \label{op_ac}\\
   z_j(t_{*}) &= z_j^{t_*},\quad j \in I_{\text{fixed}},\quad t_*\in\mathbb{R}. \label{op_rpv} 
 \end{align}
\end{subequations}
The choice of the objective criterion $\Phi$ in (\ref{op_of}) has been discussed
in \cite{Lebiedz2004c,Lebiedz2009,Lebiedz2011,Reinhardt2008}, 
\begin{equation}\label{objcrit}
 \Phi\left(z(t)\right)\coloneqq\|J_{S}\left(z(t)\right)\cdot S\left(z(t)\right)\|_2^2
\end{equation}
has been widely used recently.

Solving the optimization problem (\ref{op}), the resultant POI
$z(t_*)$ is supposed to be a good approximation to a SIM. 
In \cite{Lebiedz2011a} two different modes are presented---the {\em forward mode} ($t_*=t_0$) and the {\em reverse mode} ($t_*=t_\f$).
Both modes can be regarded as special cases of the general formulation
(\ref{op}). For the {\em reverse mode} it is shown analytically in \cite{Lebiedz2011a}, that the POI identifies the SIM exactly for an 
infinite time horizon, that is for $t_0\to-\infty$, applied to the linear
two--dimensional system (\ref{linsys}) and the nonlinear two--dimensional
Davis--Skodje test model (\ref{DS}). Accordingly, the SIM approximation error decreases
exponentially with increasing time interval ($t_\f$ fixed) which is also confirmed by numerical results for realistic chemical combustion processes including thermochemistry (see \cite{Lebiedz2011}).

Similar to the stretching--based diagnostics approach by Adrover described in
Chapter (\ref{adrover}) the optimization criterion (\ref{objcrit}) $J_{S}\left(z(t)\right)\cdot S\left(z(t)\right)$ contains the second
derivative of the state vector $J_{S}\left(z(t)\right)\cdot
S\left(z(t)\right)=\partial_t^2z(t)$. So, the
derivative--of--the--state--vector--concept is also used here.

\subsection{Boundary value view} \label{bvp}

A second fundamental concept for model reduction in kinetics with spectral gap presented in this work is the boundary--value--concept which exploits the property of attractivity of SIMs.
Provided that a SIM is globally attractive, every trajectory approaches it on infinite time horizon.
In dissipative systems assuming 
\begin{align}
 d\left(z(t_0),\text{SIM}\right)>d\left(z(t_*),\text{SIM}\right)
\end{align}
with $t_0<t_*$, $d\left(\cdot,\cdot\right)\in C^{\infty}\left(\mathbb{R}^n\times\mathbb{R}^n,\mathbb{R}\right)$, the distance function, and  
$z(t_*)=z\left(t_*-t_0,z(t_0)\right)$ (i.e.\ the solution of the initial value problem 
$\partial_tz(t)=S\left(z(t)\right),\ z(t_0)=z^{t_0}$ evaluated after a time period of $t_*-t_0$), 
the POI identifies a SIM exactly for $t_*-t_0 \rightarrow \infty$
and $d\left(z(t_0),\text{SIM}\right)=c\in\mathbb{R}$:
\begin{align}
 d\left(z(t_*),\text{SIM}\right)=0.
\end{align}
Having this in mind, the following general formulation of a boundary value problem for SIM computation
is valid
\begin{subequations}\label{BVP}
 \begin{align}
  \partial_t z(t)&=S\left(z(t)\right)\\
  z_j(t_*)&=z_j^{t_*},\quad j\in I_{\text{fixed}},\quad t_*\in\mathbb{R}\label{bovaprb}\\
  z_j(t_0)&=K_j,\quad j\notin I_{\text{fixed}},\quad K_j\in\mathbb{R}
\label{bovaprc}
 \end{align}
\end{subequations}
with $t_0<t_*$ in the {\em reverse mode} (if $t_*=t_0$ is chosen, we talk about a local method).
The crucial issue in (\ref{op}) and (\ref{zdp}) is to decide how to choose the
constant $K=(K_j)_{j\notin I_{\text{fixed}}}$. For globally attractive SIMs the choice of $K$ is without significance to obtain $\lim_{t_0\to-\infty}z_{\text{nrpv}}(t_\f)=z_{\text{rpv}}^{t_\f}$.
In contrast, in realistic chemical models the choice of $K$ plays a significant role because of additional constraints
restricting the domain where the ODE model is defined.

The conceptual idea of using such type boundary value approach for slow manifold computation is also found in \cite{Desroches2012,Guckenheimer2009,Ren2005,Ren2006a}.\\

\textbf{Boundary value approach applied to a linear system...}

\begin{figure}[h!]
  \centering
  \begin{center}
    \includegraphics[width=10cm]{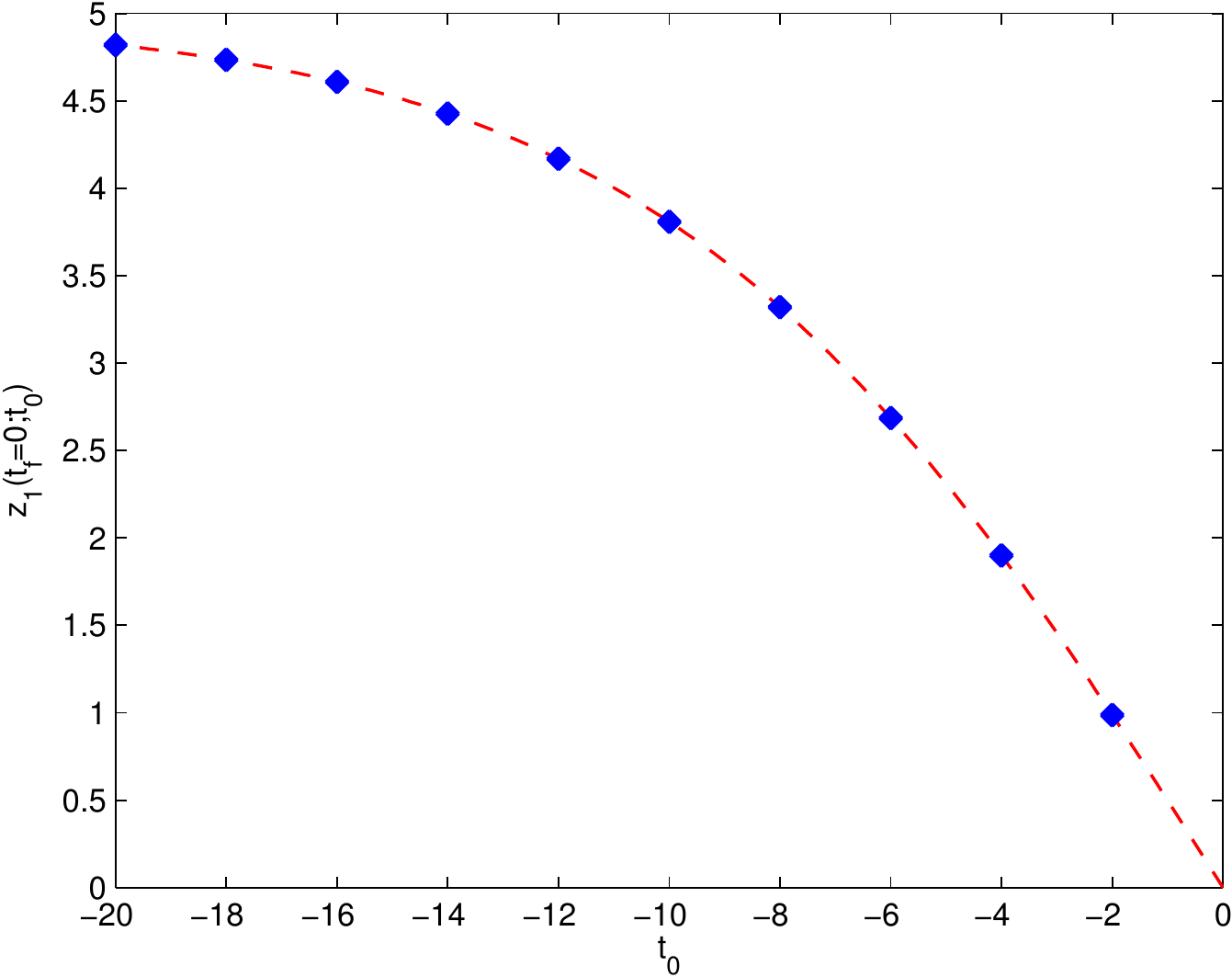}
  \end{center}
  \caption{\label{Bildchen-crop} Solutions $z_1(t_\f=0;t_0)$ of the boundary
    value formulation (\ref{BVP}) applied to the linear model (\ref{linsys}) with $\gamma=0.2$ visualized
by the blue rhombi in comparison with the analytical error (\ref{error}) (red
dashed line) as a function of $t_0$. }
\end{figure}

\begin{figure}[h!]
  \centering
  \begin{center}
    \includegraphics[width=10cm]{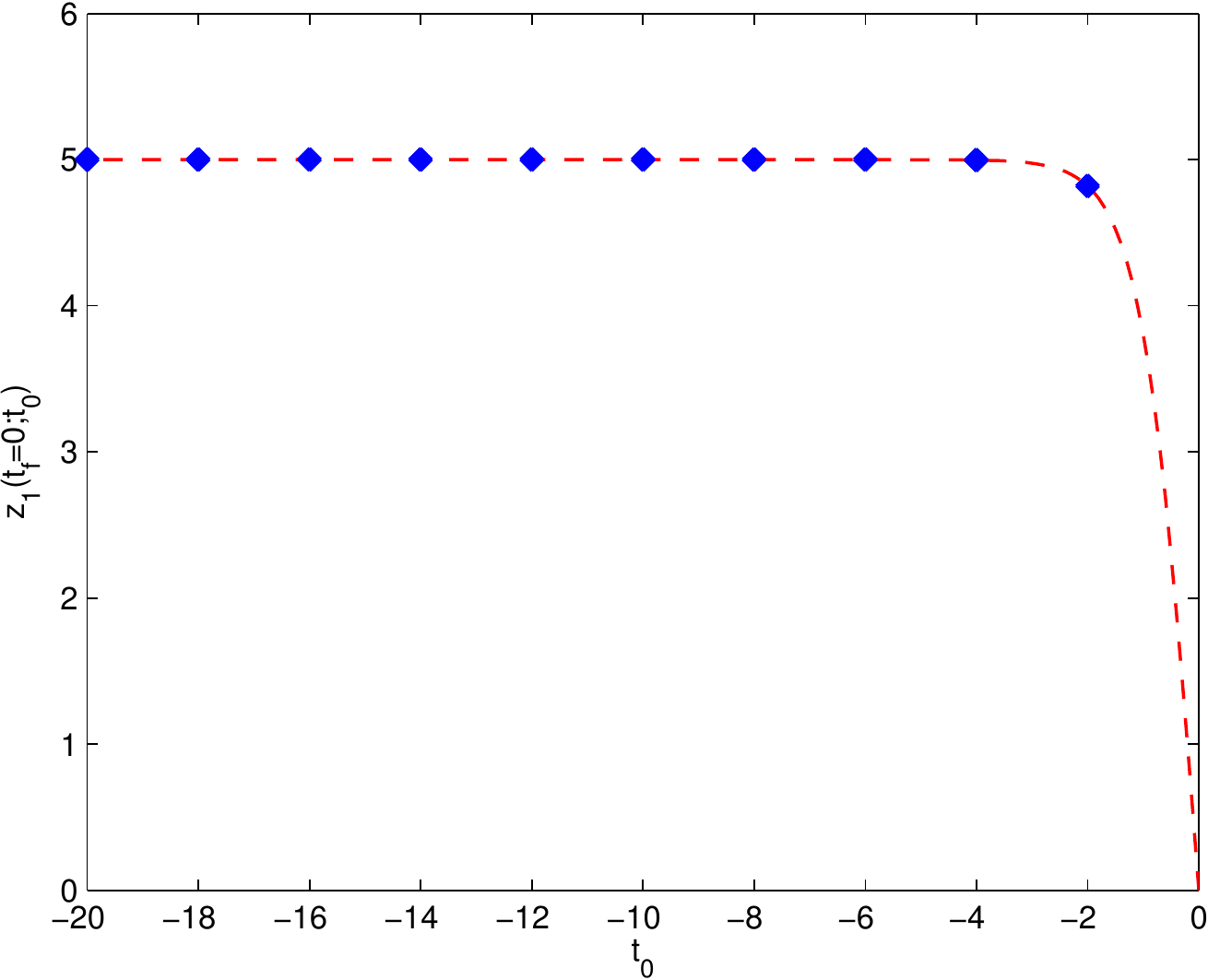}
  \end{center}
  \caption{\label{Bildchen2-crop} Solutions $z_1(t_\f=0;t_0)$ of the boundary
    value problem (\ref{BVP}) applied to the linear model (\ref{linsys}) with
    $\gamma=2.0$ visualized by the blue rhombi in comparison with the
    analytical error (\ref{error}) (red dashed line) as a function of $t_0$. }
\end{figure}

$\quad$\\\textbf{...Analytically:}\label{lin}
Again the two-dimensional linear system (\ref{linsys}) is considered, where the SIM is given by the 
first bisectrix $z_1\equiv z_2$ and the analytical solution by
\begin{subequations}
 \begin{align}
  z_1(t) &= c_1\e^{-t}+c_2\e^{(-1-\gamma) t},\quad c_1,c_2\in\mathbb{R}\\
  z_2(t) &= c_1\e^{-t}-c_2\e^{(-1-\gamma) t}.
 \end{align}
\end{subequations}
Equation (\ref{bovaprb}) yields (using $t_*=t_\f$)
\begin{align}
 z_2(t_\f)=z_2^{t_\f}=c_1\e^{-t_\f}-c_2\e^{(-1-\gamma) t_\f}
\end{align}
where $c_1(c_2)$ can be computed as
\begin{align}
c_1(c_2)=z_2^{t_\f}\e^{t_\f}+c_2\e^{-\gamma t_\f}.
\end{align}
Substituting this into Equation (\ref{bovaprc}) and using (\ref{anasolnrpv}) results in
\begin{align}
 z_1(t_0)=(z_2^{t_\f}\e^{t_\f}+c_2\e^{-\gamma t_\f})\e^{-t_0}+c_2\e^{(-1-\gamma) t_0}=K
\end{align}
and $c_2$ is obtained:
\begin{align}
c_2=\frac{K-z_2^{t_\f}\e^{t_\f}\e^{-t_0}}{\e^{-\gamma t_\f}\e^{-t_0}+\e^{(-1-\gamma)t_0}}.
\end{align}
Thus, the free variables of the POI can be computed as
\begin{equation}\label{error}
\begin{aligned}
 z_1(t_\f)&=c_1\e^{-t_\f}+c_2\e^{(-1-\gamma) t_\f}\\
\\&=z_2^{t_\f}\left[1+\frac{2\tfrac{K}{z_2^{t_\f}}\e^{(-2-\gamma)t_\f}\e^{t_0}-2\e^{(-1-\gamma)t_\f}}{\e^{(-1-\gamma)t_\f}+\e^{-t_\f}\e^{-\gamma t_0}}\right].
\end{aligned}
\end{equation}
It follows that 
\begin{subequations}\label{limes}
\begin{align} 
 \lim\limits_{\gamma\to\infty}z_1(t_\f)&=z_2^{t_\f}\\
 \lim\limits_{t_0\to-\infty}z_1(t_\f)&=z_2^{t_\f}
\end{align}
\end{subequations}
holds for every $K\in\mathbb{R}$ and $t_0<t_\f$.
This confirms the proposition, that for globally attractive systems the choice
of $K$ is not important as long as an infinite time
horizon is chosen in the {\em reverse mode} formulation.\\\\
 {\bf...Numerically:} 
Numerical experiments have been performed using \texttt{bvp4c}, a boundary value problem solver for ODEs used in MATLAB$^\circledR$.
Problem (\ref{BVP}) is implemented for the linear model (\ref{linsys}) with
$z_2(t_*)=z_2^{t_*}=5.0$, 
$t_*=t_\f=0.0$, $K=0.0$, and $t_0$ varies between $-2.0$ and $-20.0$---arbitrarily chosen values.
For Figure \ref{Bildchen-crop}, $\gamma=0.2$ is chosen and the blue rhombi
show the free variable $z_1(t_\f=0)$ resulting from the numerical solution of
the boundary value problem corresponding to the different values of $t_0$. With decreasing $t_0$, $z_1(0)$ converges to $z_1(0)=z_2(0)=5.0$ meaning that the SIM approximation improves. The red dashed line visualizes the analytical error from Equation (\ref{error}).
In Figure \ref{Bildchen2-crop} the same results for $\gamma=2.0$ are visualized.
\\

\textbf{Boundary value approach applied to the nonlinear Davis--Skodje test problem...}

\begin{figure}[h!]
  \centering
  \begin{center}
    \includegraphics[width=10cm]{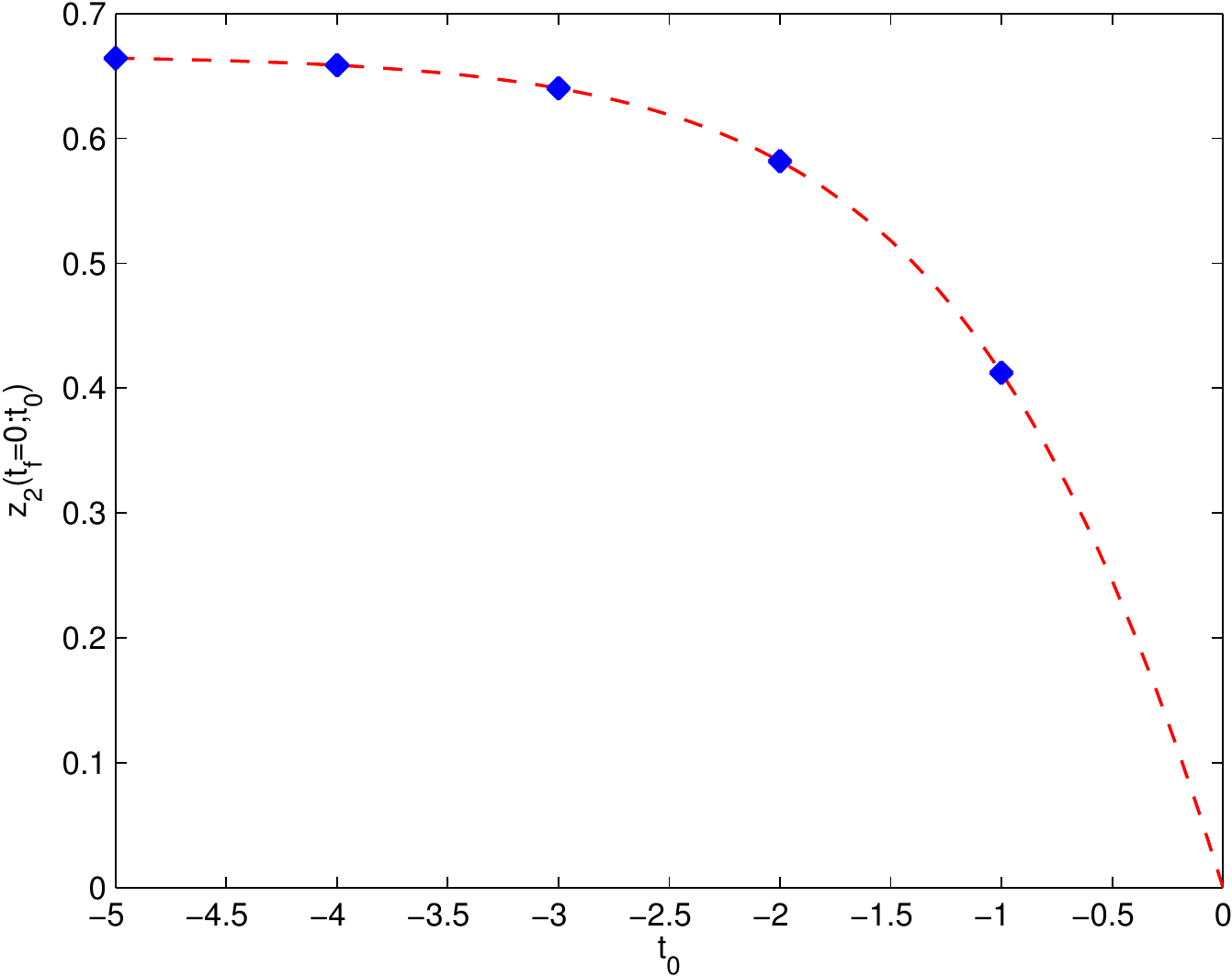}
  \end{center}
  \caption {\label{Bildchen3-crop} Solutions $z_2(t_\f=0;t_0)$ of the boundary
    value formulation (\ref{BVP}) applied to the linear model (\ref{linsys}) with $\gamma=1.2$ visualized
by the blue rhombi in comparison with the analytical error (\ref{error}) (red
dashed line) as a function of $t_0$. }
\end{figure}

\begin{figure}[h!]
  \centering
  \begin{center}
    \includegraphics[width=10cm]{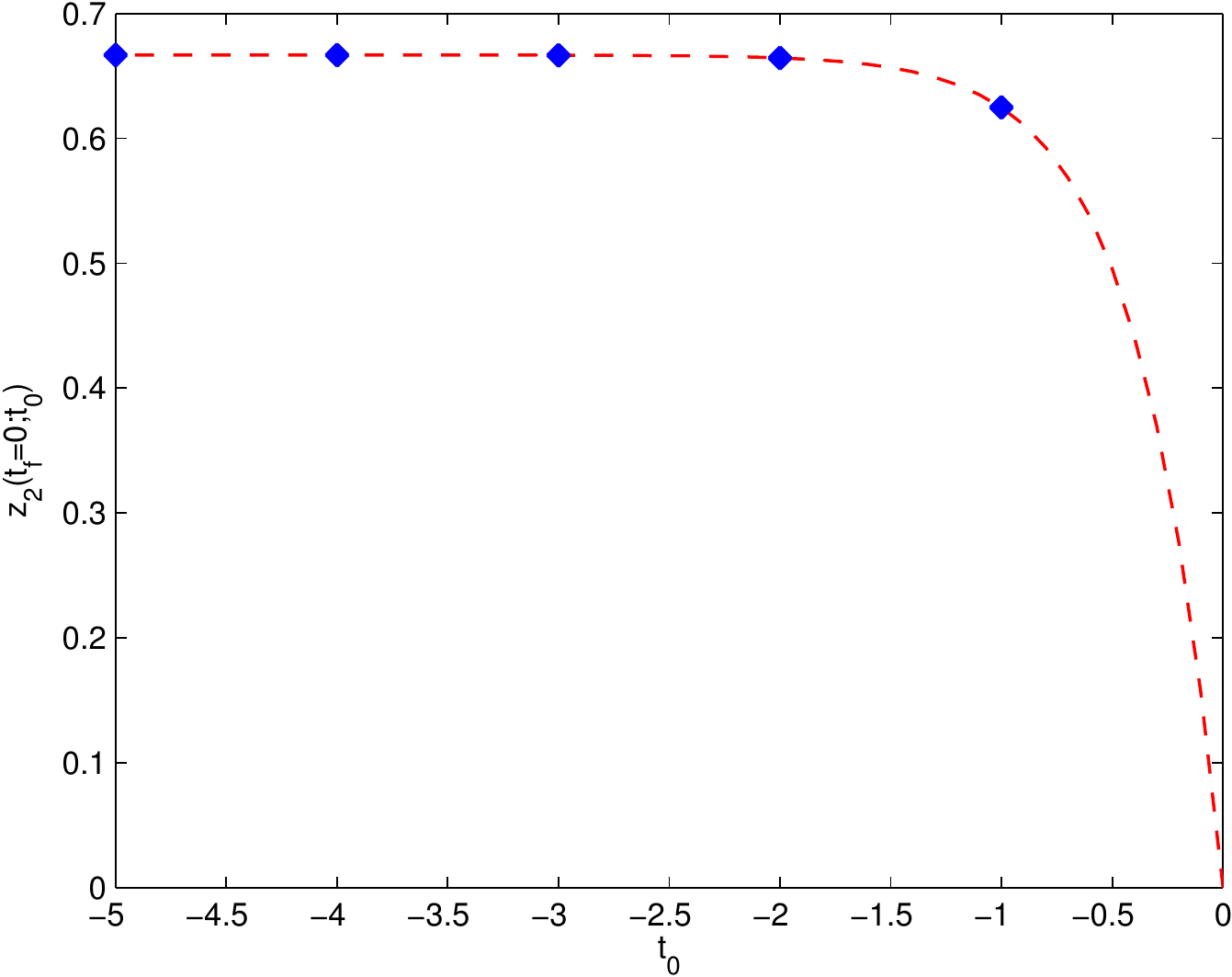}
  \end{center}
  \caption{\label{Bildchen4-crop} Solutions $z_2(t_\f=0;t_0)$ of the boundary
    value formulation (\ref{BVP}) applied to the linear model (\ref{linsys}) with $\gamma=3.0$ visualized
by the blue rhombi in comparison with the analytical error (\ref{error}) (red
dashed line) as a function of $t_0$. }
\end{figure}

$\quad$\\\textbf{...Analytically:}
The same procedure as in Section \ref{lin} is applied to the nonlinear Davis--Skodje test problem (\ref{DS}) 
where the analytically calculated SIM is given by 
\begin{align}
 z_2(t)=\frac{z_1(t)}{z_1(t)+1}
\end{align}
and the solution of the ODE system by (\ref{ansolDS}).
With $z_1$ being the RPV, the boundary values
\begin{subequations} 
 \begin{align}
  z_1(t_\f) &= z_1^{t_\f} \label{2}\\
  z_2(t_0) &= K\label{3}
 \end{align}
\end{subequations}
complete (\ref{DS}) to achieve the boundary value problem that has to be solved.
The expressions (\ref{1}) and (\ref{2}) result in
\begin{align}
 c_1=z_1^{t_\f}\e^{t_\f}
\end{align}
which, substituted into (\ref{3}), yields (using (\ref{4})):
\begin{align}
 c_2=K\e^{\gamma t_0}-\frac{z_1^{t_\f}\e^{\gamma t_0}}{z_1^{t_\f}+\e^{t_0}\e^{-t_\f}}.
\end{align}
The POI is
\begin{align}\label{errorDS}
 \begin{pmatrix}
  z_1(t_\f)\\\\z_2(t_\f)
 \end{pmatrix}
=
\begin{pmatrix}
 z_1^{t_\f}\\\\\frac{z_1^{t_\f}}{z_1^{t_\f}+1}+K\e^{\gamma t_0}\e^{-\gamma t_\f}-\frac{z_1^{t_\f}\e^{\gamma t_0}\e^{-\gamma t_\f}}{z_1^{t_\f}+\e^{t_0}\e^{-t_\f}}
\end{pmatrix}
\end{align}
which also implies that
\begin{subequations}\label{limesDS}
\begin{align} 
 \lim\limits_{\gamma\to\infty}z_{2}(t_\f)&=\frac{z_1^{t_\f}}{1+z_1^{t_\f}}\\
 \lim\limits_{t_0\to-\infty}z_{2}(t_\f)&=\frac{z_1^{t_\f}}{1+z_1^{t_\f}}
\end{align}
\end{subequations}
holds for every $K\in\mathbb{R}$ and $t_0<t_\f$.
\\\\
 {\bf...Numerically:} 
The same numerical experiment as in the linear model case is applied to the nonlinear Davis--Skodje test problem.
Here, the RPV is chosen as $z_1(t_*)=z_1(t_\f)=z_1(0)=z_1^0=2.0$ and $\gamma=1.2$ in Figure
\ref{Bildchen3-crop} and $\gamma=3.0$ in Figure \ref{Bildchen4-crop} is chosen. The constant $K$ is set to $0.0$ again and $t_0$ varies
between $-1.0$ and $-5.0$. With the analytical SIM $z_2\equiv \tfrac{z_1}{1+z_1}$, the POI should
converge towards the SIM point 
\begin{align}
 \begin{pmatrix}
  z_1^0\\\\\frac{z_1^0}{1+z_1^0}
 \end{pmatrix}
=
\begin{pmatrix}
 2\\\\\frac{2}{3}
\end{pmatrix}.
\end{align}

\subsubsection{Saddle Point Method}

Another model reduction approach exploiting a boundary value problem
formulation is the saddle point method first described by Davis and Skodje in \cite{Davis1999}.
Here, 1-D SIMs are approximated via computation and connection of fixed points
located both in physical and unphysical regions (e.g. `fixed points at
infinity') via heteroclinic orbits. This requires the use of projective
geometry with coordinate transformation 
\begin{align*}
 u_i&=\frac{z_i}{\sqrt{1+|z|^2}},\quad i=1,\dots,n\\
u_{n+1}&=\frac{1}{\sqrt{1+|z|^2}}
\end{align*}
from Euclidean space to the hyperbolic one. Here, infinity is
$u_{n+1}=0.0$. In the Davis--Skodje model, there are five fixed points, one
finite one, the equilibrium point at $(u_1,u_2)=(0,0)$ and four fixed points
at infinity $(u_1,u_2)=(0,-1)$, $(u_1,u_2)=(-1,0)$, $(u_1,u_2)=(0,1)$, and
$(u_1,u_2)=(1,0)$. By identifying the unstable manifold of a saddle point
$(u_1,u_2)=(1,0)$ the SIM is obtained by following its orbits to the stable
equilibrium point $(u_1,u_2)=(0,0)$.

This method serves as the basis for the approach developed by Al-Khateeb et al.\ \cite{Al-Khateeb2009}
where a one-dimensional SIM is defined as heteroclinic orbit---a trajectory that connects two critical points---that is locally attractive along the complete trajectory.
In \cite{Mengers2013}, this concept is enhanced to the computation of one-dimensional slow invariant manifolds of reactive systems including microscale diffusion effects.  

\subsubsection{Invariant Constrained Equilibrium Edge Preimage Curves (ICE-PIC)}

Ren et al.\ introduced the ICE-PIC approach for SIM computation in 2006 \cite{Ren2006a}. This method is based on an ICE manifold which is the union of all reaction trajectories
emanating from points located on the edge of a constrained equilibrium manifold. 
As the ICE-PIC manifold is constructed from reaction trajectories emanating
from the latter, it is invariant.
Based on this invariant constrained equilibrium edge manifold a species reconstruction can be 
done locally which means without having to generate the whole manifold in advance.
Thus, obviously the ICE-PIC method is another representative of model reduction
approaches using a boundary value problem.

\subsubsection{Trajectory Based Optimization Approach: {\em Reverse Mode}}

The {\em reverse mode} of the trajectory based optimization approach can be
formualted as 
\begin{subequations} \label{tboarm}
 \begin{equation}  \min_{z(t)} \int_{t_{0}}^{t_{\rm{f}}}\!\|J_{S}\left(z(t)\right)\cdot S\left(z(t)\right)\|_2^2 \; \textrm{d} t,\quad t_0,t_\f\in\mathbb{R},\quad t_0<t_\f 
 \end{equation}
 \textrm{subject to}
 \begin{align}
   \partial_t z(t) &= S\left(z(t)\right) \\
    z_j(t_{\rm{f}}) &= z_j^{t_{\rm{f}}},\quad j \in I_{\text{fixed}}, 
 \end{align}
\end{subequations}
where the function $g$ is omitted because of simplicity reasons.
As mentioned in Chapter \ref{tboa}, this formulation identifies SIMs exactly for an infinite time horizon, i.e.\ for $t_0\to-\infty$.
Problem (\ref{tboarm}) is (required that $t_*=t_{\rm{f}}$) a special case of
the boundary value problem (\ref{BVP}), where the objective funtional to be
minimized implicitly determines the choice of $K_j$. The variational problem
formulation in the {\em reverse mode} trajectory--based optimization approach
obviously combines ideas from both previous concepts, the state vector
derivate AND a boundary value problem.

\section{Two concepts---One approach}

A generalized ansatz to combine both concepts for model reduction presented in
the previous paragraphs is the use of derivative information in the
trajectory--based optimization method (cf.\ (\ref{tboa})). More precisely, the
general optimization problem (with $t_*=t_\f$) can be generalized to
\begin{subequations} \label{op2}
 \begin{equation}  \min_{z(t)} \int_{t_{0}}^{t_{\rm{f}}}\!\|\partial_t^mz(t)\|_2^2 \; \textrm{d} t,\quad t_0<t_\f,\quad m\in\mathbb{N} \label{op2_of}
 \end{equation}
 \textrm{subject to}
 \begin{align}
   \partial_t z(t) &= S\left(z(t)\right) \label{op2_dyn}\\
   0 &= g\left(z(t_{\f})\right) \label{op2_ac}\\
   z_j(t_{\f}) &= z_j^{t_\f},\quad j \in I_{\text{fixed}},\quad t_\f\in\mathbb{R}, \label{op2_rpv} 
 \end{align}
\end{subequations}
(note that for $m=2$ the objective criterion is equivalent to (\ref{objcrit})) and---by applying the linear model (\ref{linsys})---
the exact identification of the SIM holds for $m\to\infty$ as
it can be seen after theoretical analysis by regarding the resulting POI (compare Section 3.1 in \cite{Lebiedz2011a})
\begin{align}\label{abc}
\begin{pmatrix}
 z_1(t_\f)\\\\
 z_2(t_\f)
\end{pmatrix}
=
\begin{pmatrix}
  z_2^{t_\f}\left[1+\frac{2\e^{-2\gamma t_\f}\e^{-2t_\f}-2\e^{-2\gamma t_\f}\e^{-2t_0}}{\e^{-2\gamma t_\f}\e^{-2t_0}-\e^{-2\gamma t_\f}\e^{-2t_\f}-\xi\e^{(-1-\gamma)2t_0}+\xi\e^{(-1-\gamma)2t_\f}}\right]
\\\\z_2^{t_\f}
\end{pmatrix},\quad \xi=(-1-\gamma)^{2m-1}.
\end{align}
Thus, it holds
\begin{subequations}
\begin{align}
 \lim\limits_{\gamma\to\infty}z_1(t_\f)&=z_2^{t_\f}\\
 \lim\limits_{t_0\to-\infty}z_1(t_\f)&=z_2^{t_\f} \label{lim1}\\
 \lim\limits_{m\to\infty}z_1(t_\f)&=z_2^{t_\f}. \label{limm1}
\end{align}
\end{subequations}

Vice versa, the ZDP representing the derivative--of--the--state--vector--concept can be modified to a method using non-local trajectory information via the boundary--value--view idea.
This results in the following formulation
\begin{subequations} \label{zdp2}
 \begin{equation} \partial_t^mz_j(t)\Big|_{t=t_0}=0,\quad j\notin I_{\text{fixed}},\quad m\geq1 \label{zdp2_of}
 \end{equation}
 \textrm{subject to}
 \begin{align}
   \partial_t z(t) &= S\left(z(t)\right) \label{zdp2_dyn}\\
   0 &= g\left(z(t_{\f})\right),\quad t_\f\in\mathbb{R} \label{zdp2_ac}\\
   z_j(t_{\f}) &= z_j^{t_\f},\quad j \in I_{\text{fixed}}, \label{zdp2_rpv} 
 \end{align}
\end{subequations}
where $t_0<t_\f$ has to be fulfilled.
Theoretical analysis for the linear model (\ref{linsys}) yields
\begin{align}
 \text{POI}=
\begin{pmatrix}
 z_1(t_\f)\\\\
 z_2(t_\f)
\end{pmatrix}
=
\begin{pmatrix}
  z_2^{t_\f}\left[1+\frac{2(-1)^{m+1}\e^{(-1-\gamma)t_\f}}{(-1)^m\e^{(-1-\gamma)t_\f}+(-1-\gamma)^m\e^{-\gamma t_0}\e^{-t_\f}}\right]
\\\\z_2^{t_\f}
\end{pmatrix}
\end{align}
which implies again
\begin{subequations}\label{lim}
\begin{align} 
 \lim\limits_{\gamma\to\infty}z_1(t_\f)&=z_2^{t_\f}\\
 \lim\limits_{t_0\to-\infty}z_1(t_\f)&=z_2^{t_\f} \label{lim2}\\
 \lim\limits_{m\to\infty}z_1(t_\f)&=z_2^{t_\f}.\label{limm2}
\end{align}
\end{subequations}

In conclusion, the following formulation results as `best working', due to the fact that both the derivative--of--the--state--vector--concept and the boundary--value--view--concept are combined.
\begin{subequations} \label{op3}
 \begin{equation}  \min_{z(t)} \|\partial_t^mz(t)\|_2^2\Big|_{t=t_0},\quad m\in\mathbb{N} \label{op3_of}
 \end{equation}
 \textrm{subject to}
 \begin{align}
   \partial_t z(t) &= S\left(z(t)\right) \label{op3_dyn}\\
   0 &= g\left(z(t_{\f})\right),\quad t_\f\in\mathbb{R} \label{op3_ac}\\
   z_j(t_{\f}) &= z_j^{t_\f},\quad j \in I_{\text{fixed}}. \label{op3_rpv} 
 \end{align}
\end{subequations}

In numerical implementations for realistic chemical models, the kinetic ODE
model is only defined on a polyhedron in configuration space due to additional
constraints (e.g. species positivity, elemental mass conservation, isenthalpic
conditions) entering the optimization problem such that $t_0$ cannot be chosen
arbitrarily small.
Thus, for a good SIM approximation in realistic models the focus is on two issues to be handled:
\begin{itemize}
 \item choosing $m$ as large as practically possible (numerical computation of
   $m$-th order derivatives required),
 \item choosing $t_0$ as small as possible (with respect to the physically
   feasible domain).
\end{itemize}
The latter issue is discussed in the next section.

\subsection{Choosing $t_0$ as small as possible}

As mentioned before, the accuracy of SIM approximation in the {\em reverse mode} formulation (representing the boundary--value--view) improves with
decreasing $t_0$. In realistic chemical models the problem occurs that $t_0$ cannot get 
arbitrarily small because of additional physical constraints entering the optimization problem and 
restricting the domain where the kinetic model is defined to a polyhedron in phase space. These 
additional constraints are for instance positivity of chemical species concentrations and chemical 
element mass conservation relations. Thus, the aim is a feasible minimal choice of $t_0$, which is discussed in the following.

\begin{figure}[h!]
   \centering
    \includegraphics[width=10cm]{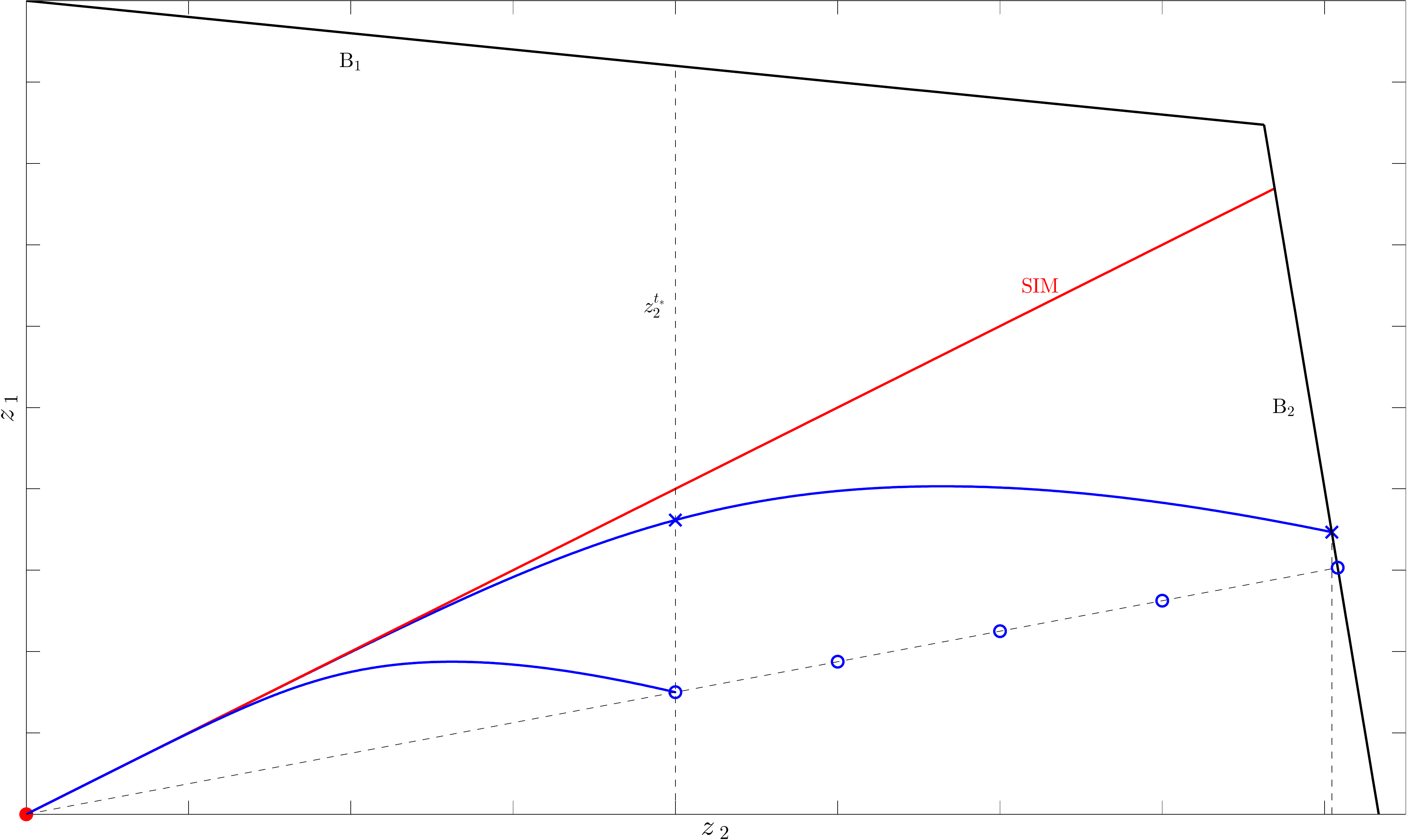}
   \caption{\label{t0-crop} Visual demonstration why a {\em reverse mode} formulation works more accurately than the corresponding local method. A polyhedron restricts the feasible area with the consequence that $t_0$ can only be chosen as small as possible within the feasibility constraints.}
\end{figure}

Figure \ref{t0-crop} exemplarily visualizes a scenario with two species ($z=(z_1,z_2)^\top$), where the phase
space polyhedron is bounded by the $z_1$-- and the $z_2$--axis ($z_1=0$ and $z_2=0$) and by two
straight lines denoted by $\text{B}_1$ and $\text{B}_2$ here. The red line refers to the SIM with
chemical equilibrium visualized by the red dot, whereas the blue lines are trajectories starting from 
specified initial values. The vertical dashed black line represents the value where the RPV $z_2$ is 
fixed at time $t=t_*$ ($\rightsquigarrow z_2^{t_*}$) and the blue circles are
the solutions of local SIM computation approaches (\ref{op3}) ($t_*=t_0=t_\f$)
with different values of $z_2^{t_*}$. The idea why the {\em reverse mode}
works better than a local method is based on the evaluation of the objective
function at time $t_0(<t_*)$. Hence, the corresponding trajectory has a time
period of $|t_0-t_*|$ to converge towards the SIM before evaluating at time
$t=t_*$ and obtaining the missing value(s) of the POI. In Figure
\ref{t0-crop}, the maximal feasible time period is represented by the blue
curve between the right cross lying on $\text{B}_2$---the result of a {\em
  reverse mode} formulation with minimal $t_0$---and the cross lying on
$z_2=z_2^{t_*}$---the point where the corresponding trajectory is evaluated at
$t=t_*$. It is obvious that the POI $z(t_*)$ has been significantly improved by using a {\em reverse mode} 
formulation compared to a local method with $z_2(t_*)=z_2^{t_*}$.

The following discussion focusses on the question, how the minimal $t_0$
could be achieved. Therefore, optimization problem (\ref{op3}) with
(\ref{linsys}) as kinetic model is regarded. For reasons of simplicity
$t_\f=0$ is chosen which is no restriction at all. Solving this problem
analytically provides formulas for the integration constants from (\ref{ansol}) depending on $t_0$ 
\begin{subequations}
\begin{align}
 \hat{c}_1&=-\frac{z_2^0}{1+\xi \e^{-2\gamma t_0}},\quad \xi=(-1-\gamma)^{2m}\\
\hat{c}_2&=z_2^0-\frac{z_2^0}{1+\xi \e^{-2\gamma t_0}}
\end{align}
\end{subequations}
which are substituted into $z_1=z_1(\hat{c}_1,\hat{c}_2)$ and $z_2=z_2(\hat{c}_1,\hat{c}_2)$ for solving the following optimization problem yielding the minimal $t_0$ that is feasible
\begin{subequations} \label{t0min}
 \begin{equation}  \min  t_0 \label{t0_of}
 \end{equation}
 \textrm{subject to}
 \begin{align}
    z_1(\hat{c}_1,\hat{c}_2)\geq 0\label{t0_1}\\
    z_2(\hat{c}_1,\hat{c}_2)\geq 0\label{t0_2}\\
    z_1(\hat{c}_1,\hat{c}_2)\leq n_1z_2(\hat{c}_1,\hat{c}_2)+b_1\label{t0_3}\\
    z_1(\hat{c}_1,\hat{c}_2)\leq n_2z_2(\hat{c}_1,\hat{c}_2)+b_2.\label{t0_4}
 \end{align}
\end{subequations}
Here, (\ref{t0_1}) and (\ref{t0_2}) are the positivity constraints of the state variables and (\ref{t0_3}) and (\ref{t0_4}) represent the restrictions $\text{B}_1$ and $\text{B}_2$ in Figure \ref{t0-crop} where the constants $n_1,n_2,b_1,b_2\in\mathbb{R}$ determine the position of these staight lines
representing a part of the boundary of the polyhedron that restricts the
domain where the kinetic model is defined. Formulas (\ref{t0_1})--(\ref{t0_4}) are examples for those additional constraints that enter the model reduction approaches above as function $g$. 

As an example, Problem (\ref{t0min}) is solved using \texttt{fmincon}---a MATLAB$^\circledR$ optimization toolbox. The following values are chosen: $\gamma=1.00$, $m=2.00$, $z_2^0=3.00$, $n_1=-2.00$, $n_2=-0.25$, $b_1=122.00$, and $b_2=111.00$. As a measure for the accuracy of the POI, the ratio $r$ between the value of the free variable of the POI and the value of the free variable of the SIM ($z_1(0)=z_2^0$) is regarded. The closer this ratio $r$ is to $r=1$, the better is the POI. Subsequently, we compare the ratio of the local method of (\ref{op3}) (that is $t_\f=t_0=0$) with the {\em reverse mode} (that is $t_0<t_\f=0$) using minimal $t_0$. Obviously, the degree of improvement depends on the parameter values chosen above, but it holds that the smaller $t_0$ the larger the improvement.
Analysis for the local method (\ref{op3}) yields
\begin{align}
 \text{POI}_{\text{loc}}=\begin{pmatrix}
      z_1(0)\\z_2^0
     \end{pmatrix}=\begin{pmatrix}
2.6471\\3.0000
\end{pmatrix}
\end{align}
which gives a ratio of $r_{\text{loc}}=\frac{2.6471}{3.0000}\approx 0.8824$.
Solving (\ref{t0min}) in the {\em reverse mode} formulation yields a minimal $t_0^{\text{min}}=-2.6056$. Using
\begin{subequations}
\begin{align}
 \hat{c}_1^{\text{min}}&=-\frac{z_2^0}{1+\xi \e^{-2\gamma t_0^{\text{min}}}}\\
\hat{c}_2^{\text{min}}&=z_2^0-\frac{z_2^0}{1+\xi \e^{-2\gamma t_0^{\text{min}}}}
\end{align}
\end{subequations}and evaluating
\begin{subequations}
 \begin{align}
  z_1(0)&=\hat{c}_1^{\text{min}}+\hat{c}_2^{\text{min}}\\
z_2(0)&=z_2^0
 \end{align}
\end{subequations}
results in
\begin{align}
 \text{POI}_{t_0^{\text{min}}}=\begin{pmatrix}
                                z_1(0)\\z_2^0
                               \end{pmatrix}
=\begin{pmatrix}
  2.9980\\3.0000
 \end{pmatrix}
\end{align}
giving a ratio of $r_{t_0^{\text{min}}}=0.9993$ which is a significant improvement compared to $r_{\text{loc}}$. The position of the polyhedron determines how small $t_0$ can be chosen. For instance, changing $b_1$ from $b_1=122.0$ to $b_1=222.0$ yields a minimal $t_0$ of $t_0^{\text{min}}=-3.2047$ and a ratio of $r_{t_0^{\text{min}}}=0.9998$. In contrast, choosing $b_1=22$ results in $t_0^{\text{min}}=0.8957$ and $r_{t_0^{\text{min}}}=0.9794$.
Apparently, the degree of improvement $|r_{t_0^{\text{min}}}-r_{\text{loc}}|$ also depends on the choice of the other variables $m,\gamma,z_2^0,n_1,n_2,b_2$.

\section{Variational principle: trajectory--based optimization approach in the light of optimal boundary control}

In the light of the boundary value problem formulation, there is a different approach to the trajectory--based optimization for model reduction in its general
formulation with $t_*=t_\f$ and without additional constraints $g$:
\begin{subequations} \label{faop}
 \begin{equation}  \min_{z(t)} \int_{t_{0}}^{t_{\rm{f}}}\!\Phi\left(z(t)\right) \; \textrm{d} t,\quad t_0<t_\f 
 \end{equation}
 \textrm{subject to}
 \begin{align}
   \partial_t z(t) &= S\left(z(t)\right) \\
    z_j(t_{\f}) &= z_j^{t_\f},\quad j \in I_{\text{fixed}},\quad t_\f\in\mathbb{R}. 
 \end{align}
\end{subequations}
The missing values of the POI, $z_j(t_\f),\ j\notin I_{\text{fixed}}$
(supposed to be an appropriate SIM approximation), are determined by the solution
of optimization problem (\ref{faop}). 
These values can be interpreted as a boundary control $u(t)$ operating at time
$t=t_\f$: 
\begin{subequations}\label{boundcontrol} 
 \begin{equation}  \min_{z(t)} \int_{t_{0}}^{t_{\rm{f}}}\!\Phi\left(z(t)\right) \; \textrm{d} t,\quad t_0<t_\f 
 \end{equation}
 \textrm{subject to}
 \begin{align}
   \partial_t z(t) &= S\left(z(t)\right) \\
    z_j(t_{\f}) &= z_j^{t_\f},\quad j \in I_{\text{fixed}},\quad t_\f\in\mathbb{R}\\
    z_j(t_{\f}) &= u(t_\f),\quad j \notin I_{\text{fixed}},\quad t_\f\in\mathbb{R}. 
 \end{align}
\end{subequations}
A criterion that automatically eliminates fast
modes as discussed in Section \ref{anatrea} is
shifted into the objective functional here. 
Formulation of the Lagrangian while introducing the Lagrange multipliers $\lambda,\mu$ leads to
\begin{align}
 J\coloneqq\int_{t_0}^{t_\f}\left[\Phi\left(z(t)\right)+\lambda^\top\left(S\left(z(t)\right)-\partial_t z(t)\right)\right]\ \d t
+\begin{pmatrix}
  \mu_{\text{nrpv}}&\mu_{\text{rpv}}
 \end{pmatrix}
\begin{pmatrix}
 z_{\text{nrpv}}(t_\f)-u(t_\f)\\z_{\text{rpv}}(t_\f)-z_{\text{rpv}}^{t_\f},
\end{pmatrix}
\end{align}
where $z_{\text{rpv}}$ denotes $(z_j)_{j \in I_{\text{fixed}}}$ and $z_{\text{nrpv}}$ denotes $(z_j)_{j \notin I_{\text{fixed}}}$.
The first variation of the Lagrangian is computed as
\begin{align}
\begin{aligned}
 \delta J=&\begin{pmatrix}
                                   \mu_{\text{nrpv}}\\\mu_{\text{rpv}}
                                  \end{pmatrix}\cdot\delta z|_{t_\f}
+\begin{pmatrix}
  z_{\text{nrpv}}(t_{\f})-u(t_\f)\\z_{\text{rpv}}(t_\f)-z_{\text{rpv}}^{t_\f}
 \end{pmatrix}\cdot\delta\mu
\\&+\int_{t_0}^{t_\f}\left[\partial_zH\cdot\delta z-\lambda\cdot\delta\dot{z}+\left(S\left(z\right)-\partial_t z\right)\cdot\delta\lambda\right]\ \d t
+\begin{pmatrix}
  -\mu_{\text{nrpv}}\\0
 \end{pmatrix}\cdot\delta u|_{t_\f}
\end{aligned}
\end{align}
with
\begin{align}
 H\coloneqq\Phi\left(z(t)\right)+\lambda^\top S\left(z(t)\right)
\end{align}
defining the Hamiltonian.
Using partial integration
\begin{align}
\int_{t_0}^{t_{\f}}\lambda\cdot\delta\dot{z}\ \d t=\lambda\cdot\delta z|_{t_{\f}}-\lambda\cdot\delta z|_{t_0}-\int_{t_0}^{t_{\f}}\dot{\lambda}\cdot\delta z\ \d t
\end{align}
leads to
\begin{align}
\begin{aligned}
 \delta J=&\left(\begin{pmatrix}
                                   \mu_{\text{nrpv}}\\\mu_{\text{rpv}}
                                  \end{pmatrix}-\begin{pmatrix}
\lambda_{\text{nrpv}}\\\lambda_{\text{rpv}}\end{pmatrix}\right)\cdot\delta z|_{t_{\f}}
+\begin{pmatrix}
  z_{\text{nrpv}}(t_{\f})-u(t_{\f})\\z_{\text{rpv}}(t_{\f})-z_{\text{rpv}}^{t_{\f}}
 \end{pmatrix}\cdot\delta\mu+\lambda\cdot\delta z|_{t_0}
\\&+\int_{t_0}^{t_{\f}}\left[\left(\partial_zH+\dot{\lambda}\right)\cdot\delta z+\left(S\left(z\right)-\partial_t z\right)\cdot\delta\lambda\right]\ \d t
+\begin{pmatrix}
  -\mu_{\text{nrpv}}\\0
 \end{pmatrix}\cdot\delta u|_{t_{\f}}.
\end{aligned}
\end{align}
The necessary optimality condition $\delta J=0$ yields the following conditions describing a boundary value problem for primal and dual variables $z(t)$ and $\lambda(t)$
\begin{subequations}\label{adjBVP}
\begin{align}
 \partial_t z(t) &= S\left(z(t)\right)\\
\partial_t \lambda(t) &= -\frac{\partial H}{\partial z}\label{adjDGL}\\
z_{\text{rpv}}(t_{\f}) &= z_{\text{rpv}}^{t_{\f}}\\
\lambda_{\text{rpv}}(t_0)&=0\\
\lambda_{\text{nrpv}}(t_0)&=0\\
\lambda_{\text{nrpv}}(t_{\f})&=0
\end{align}
\end{subequations}
with the adjoint differential equation (\ref{adjDGL}). The equations
(\ref{adjBVP}) can also be obtained from (\ref{boundcontrol}) by applying
the Pontryagin principle.

Applying this variational approach to the linear system (\ref{linsys}) 
using $\Phi\left(z(t)\right)=\|\frac{\partial^mz(t)}{\partial t^m}\|_2^2$ leads to the following Hamiltonian: 
\begin{align}
\begin{aligned}
 H=\|A^mz\|_2^2+\lambda^\top Az
=&z_1^2\left(2d_m^2+1-2d_m(-1)^m\right)\\
&+z_2^2\left(2d_m^2+1-2d_m(-1)^m\right)\\
&+z_1z_2\left(4d_m\left(-d_m+(-1)^m\right)\right)\\
&+z_1\left(-\lambda_1-\tfrac{\gamma}{2}\lambda_1+\tfrac{\gamma}{2}\lambda_2\right)\\
&+z_2\left(\tfrac{\gamma}{2}\lambda_1-\lambda_2-\tfrac{\gamma}{2}\lambda_2\right)
\end{aligned}
\end{align}
with \begin{align}
 A=
  \begin{pmatrix}
   -1-\tfrac{\gamma}{2}&\tfrac{\gamma}{2}\\\tfrac{\gamma}{2}&-1-\tfrac{\gamma}{2}
  \end{pmatrix},
 \end{align}
 \begin{align}
 A^m=
  \begin{pmatrix}
   d_m&-d_m+(-1)^m\\-d_m+(-1)^m&d_m
  \end{pmatrix},
 \end{align}
and $d_m$ being a polynomial of the form 
\begin{align}
 d_m(\gamma)=(-1)^m\left(1+\frac{m}{2}\gamma+\dots+\frac{m}{2}\gamma^{m-1}+\frac{1}{2}\gamma^m\right).
\end{align}

The adjoint differential equations can now be formulated as
\begin{subequations}
\begin{align}
\begin{aligned}
   \partial_t \lambda_1&=-\frac{\partial H}{\partial z_1}\\&=(1+\frac{\gamma}{2})\lambda_1-\frac{\gamma}{2}\lambda_2-2z_1\left(2d_m^2+1-2d_m(-1)^m\right)-z_2\left(4d_m\left(-d_m+(-1)^m\right)\right)\\
\end{aligned}
\end{align}
\begin{align}
\begin{aligned}
   \partial _t \lambda_2&=-\frac{\partial H}{\partial z_2}\\&=-\tfrac{\gamma}{2}\lambda_1+(1+\tfrac{\gamma}{2})\lambda_2-2z_2\left(2d_m^2+1-2d_m(-1)^m\right)-z_1\left(4d_m\left(-d_m+(-1)^m\right)\right)
\end{aligned}
\end{align}
\end{subequations}
with analytical solutions
\begin{subequations}
 \begin{align}
  \lambda_1(t)&=c_3\e^t+c_4\e^{(1+\gamma)t}+c_1\e^{-t}+c_2\e^{(-1-\gamma)t}\left(\frac{\left(2d_m-(-1)^m\right)^2}{1+\gamma}\right)\\
 \lambda_2(t)&=c_3\e^t-c_4\e^{(1+\gamma)t}+c_1\e^{-t}-c_2\e^{(-1-\gamma)t}\left(\frac{\left(2d_m-(-1)^m\right)^2}{1+\gamma}\right).
 \end{align}
\end{subequations}
Together with (\ref{anasol}) the Hamiltonian can be computed as
\begin{align}
\label{hamiltonian}
 H=&-2c_1c_3-2c_2c_4(1+\gamma).
\end{align}
The Hamiltonian has a remarkably simple structure.

Finally, the boundary value problem (\ref{adjBVP}) can be solved analytically leading to
\begin{subequations}\label{cs}
\begin{align}
c_1=&\frac{z_2^{t_\f}\xi\left(\e^{t_\f}\e^{(-1-\gamma)2t_0}-\e^{(-1-2\gamma)t_\f}\right)}{\xi\e^{(-1-\gamma)2t_0}-\e^{(-1-\gamma)2t_\f}(\xi+1)+\e^{-2\gamma t_\f}\e^{-2t_0}}\\\notag\\
c_2=&\frac{z_2^{t_\f}\left(\e^{(-1-\gamma)t_\f}-\e^{(1-\gamma)t_\f}\e^{-2t_0}\right)}{\xi\e^{(-1-\gamma)2t_0}-\e^{(-1-\gamma)2t_\f}(\xi+1)+\e^{-2\gamma t_\f}\e^{-2t_0}}\\\notag\\
c_3=&\frac{z_2^{t_\f}\xi\left(\e^{t_\f}\e^{(-4-\gamma)t_0}+\e^{(-1-\gamma)t_\f}\e^{-2t_0}-\e^{(1+\gamma)t_\f}\e^{(-2-\gamma)2t_0}-\e^{(-1-2\gamma)t_\f}\e^{(-2+\gamma)t_0}\right)}{\left(\xi\e^{(-1-\gamma)2t_0}-(\xi+1)\e^{(-1-\gamma)2t_\f}+\e^{-2\gamma t_\f}\e^{-2t_0}\right)\left(\e^{\gamma t_\f}-\e^{\gamma t_0}\right)}\\\notag\\
c_4=&\frac{z_2^{t_\f}\xi\left(\e^{t_\f}\e^{(-2-\gamma)2t_0}-\e^{-t_\f}\e^{(-1-\gamma)2t_0}+\e^{(-1-\gamma)t_\f}\e^{(-2-\gamma)t_0}-\e^{(1-\gamma)t_\f}\e^{(-4-\gamma)t_0}\right)}{\left(\xi\e^{(-1-\gamma)2t_0}-\e^{(-1-\gamma)2t_\f}(\xi+1)+\e^{-2\gamma t_\f}\e^{-2t_0}\right)\left(\e^{\gamma t_\f}-\e^{\gamma t_0}\right)}
\end{align}
\end{subequations}
with $\xi\coloneqq \tfrac{\left(2d_m-(-1)^m\right)^2}{1+\gamma}$ and $I_\text{fixed}=\{2\}$.
The missing value of the POI $z_1(t_\f)$ can now be computed as
\begin{align}
 z_1(t_\f)=
z_2^{t_{\f}}\left[1+\underbrace{\left(\frac{2\e^{\left(-1-\gamma\right)2t_{\f}}-2\e^{-2\gamma t_{\f}}\e^{-2t_0}}{\e^{-2\gamma t_{\f}}\e^{-2t_0}+\xi\e^{\left(-1-\gamma\right)2t_0}-\left(\xi+1\right)\e^{\left(-1-\gamma\right)2t_{\f}}}\right)}_{\eqqcolon \chi}\right]
\end{align}
where the error term $\xi$ is equivalent to (\ref{abc})
where the POI is computed by directly solving the optimization problem (\ref{faop}) analytically. 
Consequently, it holds
\begin{subequations}
\begin{align} 
 \lim\limits_{\gamma\to\infty}z_1(t_\f)&=z_2^{t_\f}\\
 \lim\limits_{t_0\to-\infty}z_1(t_\f)&=z_2^{t_\f}.
\end{align}
\end{subequations}

The boundary control formulation could be exploited for efficient numerical
implementation of trajectory-based slow manifold computation since the dual
variable $\lambda$ can be used to compute the gradient of the objective
function with respect to the system state and thus yields derivative
information by a single numerical integration of the adjoint differential
equation (see \cite{Cao2003}, Chapter 2.1).

\section{Hamilton's principle, (partial) integrability and symmetry issues in
  the search for an exact objective functional}

Based on empirical work of Lebiedz and Reinhardt \cite{Reinhardt2008a, Reinhardt2008}
and their results concerning the use of an additive term in the objective function
(\ref{op_of}), (\ref{objcrit}), and due to the non-physical fact that 
$\lim_{t_0\to-\infty} \| H \| =\infty$ for the `energy-like' Hamiltonian $H$ 
(\ref{hamiltonian}) with $c_1-c_4$ substituted from (\ref{cs}), we conjecture
a possible lack of some additional term in the formulation (\ref{op_of}),
(\ref{objcrit}) in order to achieve an exact identification of slow manifolds
via a finite--time--horizon, finite--derivative--order variational approach without
using limit arguments. This proposition is motivated by analogy reasoning with
respect to Hamilton's principle -- the principle of stationary action -- in
classical mechanics and its conceptual generalization to disspipative systems
where the `generalized forces' cannot be derived from a potential, see e.g.\ \cite{Santilli1978,Santilli1983}. The full system information is collected in the
functional of the variational problem and encoded in a single function,
the Lagrangian $\mathcal{L}\left(z(t),\partial_tz(t),t\right)$. In classical
mechanics, the Lagrangian is characterized by the difference of kinetic and
potential energy $T\left(\partial_z(t),t\right)-V\left(z(t),t\right)$, which
in our case suggests to consider the following formulation of the objective
function 
(\ref{op_of}) 
\begin{align}\label{ks}
 \min_{z(t)} \int_{t_{0}}^{t_{\rm{f}}}\!k_1\|\partial_tz(t)\|_2^2-k_2\|z(t)\|_2^2 \; \textrm{d} t
\end{align}
with constants $k_1,k_2\in\mathbb{R}$ determining the `quality' of SIM
approximation. The first integrand term corresponds to some `generalized kinetic
energy' (proportional to squared velocity) and the second to some `generalized
potential energy' (proportional to the squared deviation of the state $z(t)$
from equilibrium $(0,0)$ here).
As mentioned before, exact SIM identification requires $c_2=0$ for the two
test models analyzed in Chapter \ref{anatrea} which can be achieved by
$k_2=1$ in the linear model and $k_2=\tfrac{\gamma}{z_1(t)+1}$ in the
Davis--Skodje test model for fixed $k_1=1$.
Moreover, for a three-dimensional linear model given by
\begin{subequations}\label{linsys3D}
\begin{align}
 \partial_tz_1(t)&=\left(-1-\frac{\gamma_1}{4}-\frac{\gamma_2}{2}\right)z_1(t)+\frac{\sqrt{2}\gamma_1}{4}z_2(t)+\left(\frac{\gamma_2}{2}-\frac{\gamma_1}{4}\right)z_3(t)\\
 \partial_tz_2(t)&=\frac{\sqrt{2}\gamma_1}{4}z_1(t)-\left(1+\frac{\gamma_1}{2}\right)z_2(t)+\frac{\sqrt{2}\gamma_1}{4}z_3(t)\\
 \partial_tz_3(t)&=\left(-\frac{\gamma_1}{4}+\frac{\gamma_2}{2}\right)z_1(t)+\frac{\sqrt{2}\gamma_1}{4}z_2(t)-\left(1+\frac{\gamma_1}{4}+\frac{\gamma_2}{2}\right)z_3(t),\quad\gamma_1,\gamma_2>0,
\end{align}
\end{subequations}
with $\gamma_1,\gamma_2\in\mathbb{R}$, $t\in\mathbb{R}$, $z_1,z_2,z_3\in C^{\infty}\left(\mathbb{R},\mathbb{R}\right)$, and analytical solutions
\begin{subequations}\label{ansol3D}
 \begin{align}
  z_1(t) &= c_1\e^{-t}+c_2\e^{(-1-\gamma_1) t}+c_3\e^{(-1-\gamma_2) t}\\
  z_2(t) &= \sqrt{2}c_1\e^{-t}-\sqrt{2}c_2\e^{(-1-\gamma_1)t}\\
  z_3(t) &= c_1\e^{-t}+c_2\e^{(-1-\gamma_1) t}-c_3\e^{(-1-\gamma_2) t},\quad c_1,c_2,c_3\in\mathbb{R},
 \end{align}
\end{subequations}
a two-dimensional SIM can be computed exactly by using (\ref{ks}) with $k_1=k_2=1$ as well.
Here, the slow manifold which is spanned by the two eigenvectors corresponding to the slow eigenvalues of system (\ref{linsys3D}) is represented by $z_2(t)=h\left(z_1(t),z_2(t)\right)=\tfrac{z_1(t)+z_3(t)}{\sqrt{2}}$.
To find a general characterization of $k_1$ and $k_2$ or a suitable general 
form of the Lagrangian (the inverse problem in the calculus of variations, see
e.g. \cite{Santilli1978,Morandi1990}) leading to an exact SIM identification
in a variational approach without using limiting arguments would be an
important issue for 
model reduction in chemical kinetics. We believe that a Hamiltonian
variational formulation might turn out to be a promising approach towards this
goal. The Hamiltonian viewpoint offers an interesting new perspective on
slow invariant manifolds. According to Pontryagin's maximum principle
\cite{Pontryagin1960}, the Hamiltonian $H$ is constant along the optimal
solution of a variational problem with given non-explicitly time-dependent
Lagrangian and non-holonomic constraints given by autonomous ordinary
differential equations. Conserved properties are closely related to the issue
of (partial) integrability and the existence of various types of first
integrals of dynamical systems.  
If an approximated slow invariant manifold can be computed as a solution of a
variational problem, it is obviously related to the existence of a
conservation relation along trajectories on the manifold. According to
Noether's theorem \cite{Noether1918} conservation relations are related to
symmetries of the Lagrangian. Symmetries are generally essential in the
macroscopic modeling of multiscale problems because non-trivial macroscopic
dynamics can only occur if `microscopic modes' do not cancel out completely,
which requires to existence of symmetries. 

\begin{figure}[htbp]
  \centering
  \subfigure[Linear model]{
    \label{Labelname 1}
    \includegraphics[width=6cm]{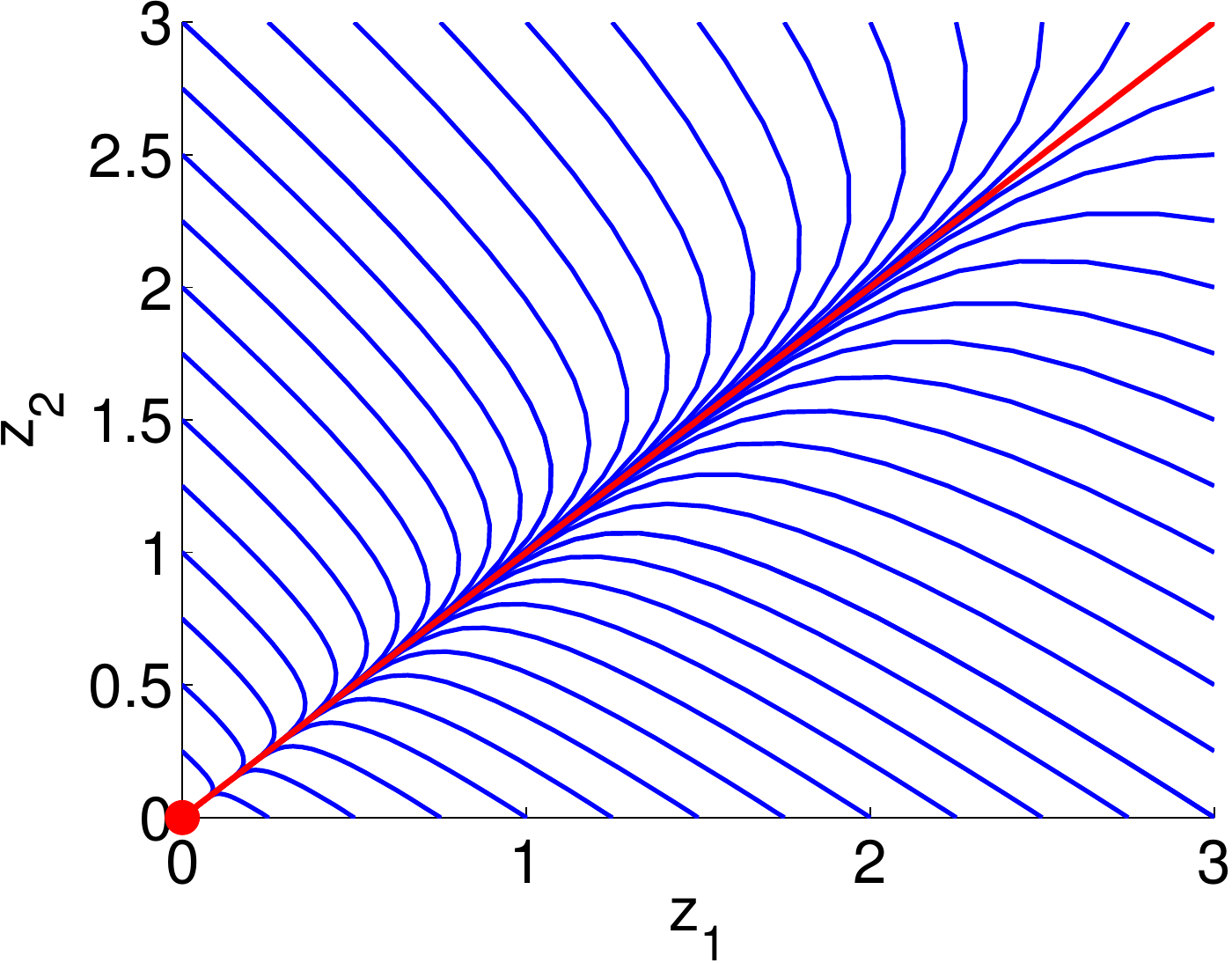} 
  }
  \subfigure[Davis--Skodje model]{
    \label{Labelname 2}
    \includegraphics[width=6cm]{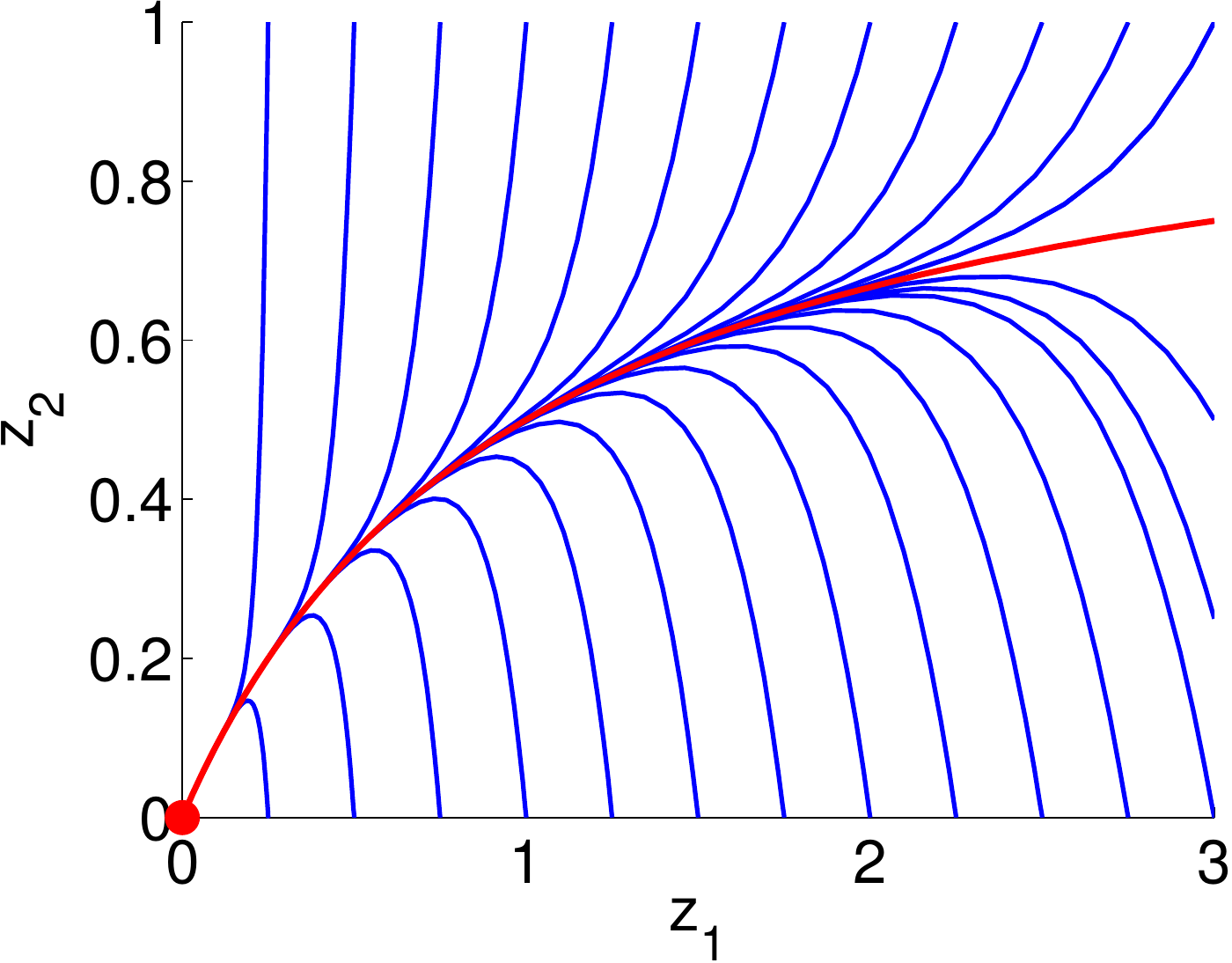} 
  }
  \caption{The red dots represent the equilibrium point, the red curves the SIM, and the blue curves are solution trajectories concerning the underlying model equations.}
  \label{LinMod&DS}
\end{figure}

In the case of the 2-D linear model
and the 2-D Davis--Skodje model analyzed on this article, it seems that
the SIM, the tangent space of the SIM respectively, correspond to the symmetry
axes of local mirror symmetry on the manifold of solution trajectories of the
ODE systems (see Fig.\ \ref{LinMod&DS}). We consider
these issues to be important for a deep understanding of the 
unifying elements of various model reduction approaches computing slow
invariant attracting manifolds in chemical kinetics.

\begin{thebibliography}{11} 

\bibitem{Adrover2007a}
 A.~Adrover, F.~Creta, S.~Cerbelli, M.~Valorani, and M.~Giona, {\em The structure of slow invariant manifolds and their bifurcational
	routes in chemical kinetic models}, 
Comput. Chem. Eng., 31 (2007), pp.~1456--1474.

\bibitem{Adrover2007}
 A.~Adrover, F.~Creta, M.~Giona, and M.~Valorani, {\em Stretching-based
  diagnostics and reduction of chemical kinetic models with diffusion}, 
J. Comput. Phys., 225 (2007), pp.~1442--1471.

\bibitem{Pontryagin1960}
 V.G.~Boltyanskii, R.V.~Gamkrelidze, and L.S.~Pontryagin, {\em Theory of optimal processes. I. The maximum principle}, 
Izv. Akad. Nauk SSSR Ser. Mat., 24 (1960), pp.~3--42.

\bibitem{Al-Khateeb2009}
 A.N.~Al-Khateeb, J.M.~Powers, S.~Paolucci, A.J.~Sommese, J.A.~Diller, J.D.~Hauenstein, and J.D.~Mengers,
{\em One--dimensional slow invariant manifolds for spatially homogenous reactive systems},
J. Chem. Phys., 131 (2009), p.~024118
 
\bibitem{Bodenstein1913}
 M. Bodenstein, {\em Eine {T}heorie der photochemische {R}eaktionsgeschwindigkeiten}, 
Z. Physik. Chem., 85 (1913), pp.~329--397.

\bibitem{Cao2003}
 Y. Cao, S. Li, L. Petzold and R. Serban, {\em Adjoint sensitivity analysis for differential--algebraic equations: The adjoint DAE system and its numerical solution}, 
SIAM J. Sci. Comput., 24 (2003), pp.~1076--1089.

\bibitem{Chapman1913}
 D.L. Chapman and L.K. Underhill, {\em The interaction of chlorine and hydrogen. {T}he influence of mass.}, 
J. Chem. Soc., 103 (1913), pp.~496--508.

\bibitem{Chiavazzo2011}
E.~Chiavazzo and I.~Karlin, {\em Adaptive simplification of complex multiscale systems.} Phys.
Rev. E, 83 (2011), p. 036706.

\bibitem{Darboux1878}
 G. Darboux, {\em Sur les équations différentielles algébriques du premier ordre et du premier degré}, 
Bull. Sci. Math., Sr. 2 (1878), pp. 60--96, pp. 123--143, pp. 151--200.

\bibitem{Davis1999}
 M.J. Davis and R.T. Skodje, {\em Geometric investigation of
  low--dimensional manifolds in systems approaching equilibrium}, 
J. Chem. Phys., 111 (1999), pp.~859--874.

\bibitem{Desroches2012}
 M. Desroches, J. Guckenheimer, B. Krauskopf, C. Kuehn, H.M. Osinga, and M. Wechselberger, {\em Mixed-mode oscillations with
multiple time scales}, 
SIAM Rev., 54 (2012), pp.~211--288.

\bibitem{Gear2005}
 C.W. Gear, T.J. Kaper, I.G. Kevrekidis, and A.~Zagaris, {\em Projecting
  to a slow manifold: Singularly perturbed systems and legacy codes}, 
SIAM J. Appl. Dyn. Syst., 4 (2005), pp.~711--732.

\bibitem{Goussis2012}
 D.A.~Goussis, {\em Quasi steady state and partial equilibrium approximations: their relation and their validity}, 
Combust. Theor. Model., 16 (2012), pp.~869--926.

\bibitem{Ginoux2008}
 J.M. Ginoux, B. Rossetto, and L.~Chua, {\em Slow Invariant Manifolds as Curvature of the Flow of Dynamical Systems}, 
Int. J. Bifurcat. Chaos, 18 (2008), pp.~3409--3430.

\bibitem{Guckenheimer2009}
 J.~Guckenheimer and C.~Kuehn, {\em Computing slow manifolds of saddle type}, 
SIAM J. Appl. Dyn. Syst., 8 (2009), pp.~854--879.

\bibitem{Jones1994}
 C.K.R.T.~Jones, {\em Geometric singular perturbation theory}, 
Lect. Notes Math., 1609 (1994), pp.~44--118.

\bibitem{Kaper1999}
 T.J.~Kaper, {\em An introduction to geometric methods and dynamical systems theory for singular perturbation problems}, 
in Analyzing Multiscale Phenomena Using Singular Perturbation Methods, Proc. Symp. Appl. Math., 56, 1999, R. E. O'Malley, Jr., and J. Cronin, eds., Am. Math. Soc., Providence, RI, pp. 85--132. 

\bibitem{Kaper2002}
 H.G.~Kaper and T.J.~Kaper, {\em Asymptotic analysis of two reduction methods for systems of chemical reactions}, 
Physica D, 165 (2002), pp.~66--93.

\bibitem{Kevrekidis2003}
 I.G.~Kevrekidis, C.W.~Gear, J.M.~Hyman, P.G.~Kevrekidis, O.~Runborg, and
C.~Theodoropoulos, {\em Equation-free, coarse-grained multiscale computation: Enabling
microscopic simulators to perform system-level analysis}, Commun. Math. Sci., 1 (2003), pp.~715--762.


\bibitem{Lam1985}
 S.H. Lam, {\em Recent Advances in the Aerospace Sciences}, Plenum Press,
  New York and London, 1985, ch.~Singular Perturbation for Stiff Equations
  using Numerical Methods, pp.~3--20.

\bibitem{Lam1994}
 S.H. Lam and D.A. Goussis, {\em The {CSP} method for simplifying
  kinetics}, 
Int. J. Chem. Kinet., 26 (1994), pp.~461--486.

\bibitem{Lebiedz2004c}
 D.~Lebiedz, {\em Computing minimal entropy production trajectories: An
  approach to model reduction in chemical kinetics}, 
J. Chem. Phys., 120 (2004), pp.~6890--6897.

\bibitem{Lebiedz2006b}
 D.~Lebiedz, V.~Reinhardt, and J.~Kammerer, {\em Novel trajectory based
  concepts for model and complexity reduction in (bio)chemical kinetics}, in
  Model reduction and coarse--graining approaches for multi-scale phenomena,
  A.~N. Gorban, N.~Kazantzis, I.~G. Kevrekidis, and C.~Theodoropoulos, eds.,
  Springer, Berlin, 2006, pp.~343--364.

\bibitem{Lebiedz2009}
 D.~Lebiedz, V.~Reinhardt, and J.~Siehr, {\em Minimal curvature
  trajectories: Riemannian geometry concepts for model reduction in chemical
  kinetics}, 
J. Comp. Phys., 229 (2010), pp.~6512--6533.

\bibitem{Lebiedz2010a}
 D.~Lebiedz, {\em Entropy--related extremum principles for model reduction of dynamical
	systems}, 
Entropy, 12 (2010), pp.~706--719.

\bibitem{Lebiedz2011}
 D.~Lebiedz, V.~Reinhardt, J.~Siehr, and J.~Unger, {\em Geometric criteria for model reduction in chemical kinetics via optimization
	of trajectories}, in Coping with Complexity: Model Reduction and Data Analysis, 
A.~N. Gorban, D. Roosepp, eds., Springer, 2011, pp.~241--252.

\bibitem{Lebiedz2011a}
 D.~Lebiedz, J.~Siehr, and J.~Unger, {\em A variational principle for computing slow invariant manifolds in
	dissipative dynamical systems}, 
SIAM J. Sci. Comput., 33
(2011), pp.~703--720.

\bibitem{Maas1992}
 U.~Maas and S.B. Pope, {\em Simplifying chemical kinetics: Intrinsic
  low--dimensional manifolds in composition space}, 
Combust. Flame, 88 (1992), pp.~239--264.

\bibitem{Mease2012}
 K.D.~Mease, U. Topcu, E. Aykutlug and M. Maggia, {\em Characterizing two-timescale nonlinear dynamics using
finite-time Lyapunov exponents and vectors}, 
arXiv:0807.0239v2 [math.DS], 2012. 

\bibitem{Michaelis1913}
 L.~Michaelis and M.L. Menten, {\em Die {K}inetik der {I}nvertinwirkung}, 
Biochem. Z., 49 (1913), pp.~333--369.

\bibitem{Morandi1990}
 G.~Moradni, C.~Ferrario, G.L.~Vecchio, G.~Marmo, and C.~Rubano, {\em
  The inverse problem in the calculus of variations and the geometry of the
  tangent bundle}, 
Phys. Rep. 188 (1990), pp.~147--284.

\bibitem{Noether1918}
 E.~Noether, {\em Invariante Variationsprobleme}, 
Nachr. D. K\"onig. Gesellsch. D. Wiss. Zu G\"ottingen, Math-phys. Klasse, (1918), pp.~235--257.

\bibitem{Mengers2013}
 J.D.~Mengers and J.M.~Powers, {\em One-Dimensional Slow Invariant Manifolds for Fully Coupled Reaction and
Micro-scale Diffusion}, 
SIAM J. Appl. Dyn. Syst., 12 (2013), pp.~560--595.

\bibitem{Reinhardt2008a}
 V.~Reinhardt, {\em On the application of trajectory--based optimization for nonlinear
	kinetic model reduction}, 
Ph.D. thesis, University of Heidelberg, Heidelberg, Germany (2008).

\bibitem{Reinhardt2008}
 V.~Reinhardt, M.~Winckler, and D. Lebiedz, {\em Approximation of slow attracting manifolds in chemical kinetics by
	trajectory--based optimization approaches}, 
J. Phys. Chem. A, 112 (2008),
  pp.~1712--1718.

\bibitem{Ren2005}
 Z.~Ren and S.B. Pope, {\em Species reconstruction using pre-image
  curves}, in Proc. Comb. Inst., 30 (2005), pp.~1293--1300.

\bibitem{Ren2006a}
 Z.~Ren, S.B. Pope, A.~Vladimirsky, and J.M. Guckenheimer, {\em The
  invariant constrained equilibrium edge preimage curve method for the
  dimension reduction of chemical kinetics}, 
J. Chem. Phys., 124 (2006), p.~114111.

\bibitem{Rossetto1986}
 B. Rossetto, {\em Trajectoires lentes des systèmes dynamiques}, in Proceedings of the 7th International Conference on Analysis and Optimization of Systems, Antibes, Juan-les-Pins, France, 
Springer, 1986, pp.~630--645.

\bibitem{Roussel2006}
 M.R.~Roussell, and T. Tang, {\em The functional equation truncation method for approximating slow
invariant manifolds: A rapid method for computing intrinsic
low-dimensional manifolds}, 
J. Chem. Phys., 125 (2006), p.~214103.

\bibitem{Roussel2012}
 M.R.~Roussell, {\em Further studies of the functional equation truncation approximation}, 
Can. Appl. Math. Q., 20 (2012), pp.~209--227.

\bibitem{Santilli1978}
 R.M.~Santilli, {\em Foundations of Theoretical Mechanics I: The inverse
  problem in Newtonian mechanics.}, Springer, New York, 1978

\bibitem{Santilli1983}
 R.M.~Santilli, {\em Foundations of Theoretical Mechanics II: Birkhoffian
  Generalization of Hamiltonian Mechanics.}, Springer, New York, 1983

\bibitem{Singh2002}
 S.~Singh, J.M. Powers, and S.~Paolucci, {\em On slow manifolds of
  chemically reactive systems}, 
J. Chem. Phys., 117 (2002), pp.~1482--1496.

\bibitem{Theodoropoulos2000}

 C.~Theodoropoulos, Y.-H.~Qian, and I.G.~Kevrekidis, {\em Coarse stability and bifurcation
analysis using time-steppers: A reaction-diffusion example},  Proc. Natl. Acad. Sci., 97 (2000), pp.~9840--9843.

\bibitem{Valorani2009}
 M.~Valorani and S.~Paolucci, {\em The g-scheme: A framework for multi-scale adaptive model
reduction}, J. Comp. Phys., 228 (2009), pp. 4665--4701.

\bibitem{Zagaris2009}
 A.~Zagaris, C.W. Gear, T.J. Kaper, and Y.G. Kevrekidis, {\em Analysis of
  the accuracy and convergence of equation--free projection to a slow manifold},
  ESAIM: Math. Model. Num., 43 (2009), pp.~757--784.
\end{thebibliography}
\end{document}